\documentclass[10pt]{article}

\usepackage{amsmath}
\usepackage{amsfonts}
\usepackage{url}
\usepackage{amsthm}
\usepackage{color}
\usepackage{graphicx}
\usepackage{tcolorbox}
\usepackage{authblk}

\theoremstyle{remark}

\newtheorem{myrem}{Remark}

\usepackage[margin=2cm]{geometry}
\title{Models for plasma kinetics during  simultaneous therapeutic plasma exchange and extracorporeal membrane oxygenation}
\author{Charles Puelz$^1$ \hspace{0.15cm} Zach Danial$^1$ \hspace{0.15cm}  Jon Marinaro$^2$ \hspace{0.15cm}  \\    Jay Raval$^3$ \hspace{0.15cm} Boyce E. Griffith$^{4,5,6,7}$ \hspace{0.15cm} Charles S. Peskin$^1$}
\date{%
    $^1$Courant Institute of Mathematical Sciences, New York University\\%
    $^2$Department of Emergency Medicine, University of New Mexico\\
    $^3$ Department of Pathology, University of New Mexico \\
    $^4$Departments of Mathematics, Applied Physical Sciences, and Biomedical Engineering, University of North Carolina\\
    $^5$Carolina Center for Interdisciplinary Applied Mathematics, University of North Carolina \\
    $^6$Computational Medicine Program, University of North Carolina \\
    $^7$McAllister Heart Institute, University of North Carolina\\[2ex]%
    \today
}
\begin{document}

\maketitle

\abstract{This paper focuses on the derivation and simulation of mathematical models describing new plasma fraction in blood for patients undergoing simultaneous extracorporeal membrane oxygenation and therapeutic plasma exchange.  Models for plasma exchange with either veno--arterial or veno--venous extracorporeal membrane oxygenation are considered.  Two classes of models are derived for each case, one in the form of an algebraic delay equation and  another in the form of a system of delay differential equations. In special cases, our models reduce to single compartment ones for plasma exchange that have been validated with experimental data   \cite{Randerson82}.  We also show that the algebraic delay equations are forward Euler discretizations of the delay differential equations, with timesteps equal to transit times through model compartments.  Numerical simulations are performed to compare different model types, to investigate the impact of plasma device port switching on the efficiency of the exchange process, and to study the sensitivity of the models to their parameters.}

\section{Introduction}
Therapeutic Plasma Exchange (TPE) refers to the removal and possible exchange of a patient's blood plasma.  This procedure treats a variety of conditions, from renal diseases to neurological disorders~\cite{Madore96, Chhibber12}. It has very recently received attention as a possible therapy for COVID--19 \cite{kesici20, Keith20, Duan20}.  TPE is typically done by connecting the patient's circulation to a device that slowly draws in blood, separates out old plasma, adds in new plasma or some other replacement fluid, and puts this mixture back into the native circulation.  Extracorporeal Membrane Oxygenation (ECMO) is a procedure used for reoxygenating blood with an external circuit.  There are two main types of ECMO, veno--venous and veno--arterial, and they are distinguished from each other by the placement locations of the inlet and return lines within the patient's circulation~\cite{Pavlushkov17}.  Veno--venous ECMO draws from and replaces blood within the venous circulation.  It is traditionally used when the patient's heart is functioning well enough to perfuse their organs.  Veno--arterial ECMO bypasses both the heart and lungs and is used in cases where the patient needs both reoxygenation and circulatory support.  In this case, blood is drawn from the veins and re--enters the native circulation through the arteries.

The clinical scenario of interest in this paper is  simultaneous ECMO and plasma exchange, which is usually done by connecting the TPE device directly to the ECMO circuit, often at the ECMO inlet line~\cite{Chong17, Dyer14, Laverdure18}.  This technique has been used to treat conditions such as multisystem organ failure or the rejection of a transplanted organ~\cite{Jhang07, Chong17, Dyer14}.  There are multiple parameters that affect the plasma exchange procedure in these cases, including inlet and return flow rates for the ECMO and TPE devices, fraction of native cardiac output supported by the ECMO circuit, and recirculation of new plasma.  This paper develops mathematical models that help quantify the impact of these parameters on plasma exchange done with ECMO.  We derive two classes of models that each describe fraction of new plasma over time, and we consider both veno--venous and veno--arterial ECMO.  The approaches used in this paper can be applied to other apheresis techniques done in conjunction with ECMO.

To our knowledge, there are no mathematical models for simulating simultaneous ECMO and plasma exchange besides our previous work~\cite{Puelz20}.  This paper provides mathematical justification for our previous model, extends this approach to both veno--arterial and veno--venous ECMO, and derives a new class of related models. There has been some research focused on separately modeling either plasma exchange or ECMO procedures.  A compartmental model of plasma exchange kinetics was described by Kellogg and Hester~\cite{Kellogg88}.  They derived equations for plasma concentrations in both the intra-- and extra--vascular spaces.  A single compartment model for plasma exchange was proposed by Randerson et al.~\cite{Randerson82}.  Our approach is similar to the models introduced in these papers in that it describes plasma kinetics through several distinguished compartments, but our focus is on the interaction between the heart/lung, peripheral, and ECMO compartments.  When these compartments are lumped together, we show that our model reduces to that of Randerson et al.~\cite{Randerson82}. Our approach is also able to describe possible recirculation of new plasma through these various spaces.  As far as  simulations for ECMO procedures, we mention the work by Zanella et al., which focuses on oxygenation in veno--venous ECMO \cite{Zanella16}.  

This paper provides derivations for the two classes of models: algebraic delay equations and systems of delay differential equations.  We also consider two possible configurations for the TPE device, either ``typical'' and ``switched.''  These configurations refer to the orientation of the TPE device ports with respect to the ECMO device flow direction.  In the typical configuration, the inlet line for the TPE device is upstream from the return line, and in the switched configuration these lines are reversed.  A comparison of these two configurations was the focus of our initial modeling effort, and we further study these two configurations in this paper \cite{Puelz20}.  We also present numerical results comparing different classes of models and different configurations of the TPE device.

\section{Mathematical model descriptions}
\label{sec:models}

Our approach is to divide the system consisting of the patient, ECMO circuit, and TPE device into three compartments: (1) the heart and lungs; (2) the peripheral organs; and (3) the ECMO circuit with the attached TPE device.  Each compartment has an associated volume and flow rate of blood. The ratio of of volume to flow determines the transit time of a parcel of new plasma through the given compartment.  In particular, fixing two of the parameters for a given compartment (transit time, flow, or volume) will determine the third.  Flows through compartments will always be specified, and either transit times or volumes will be given, depending on available data.  Models are derived in this section for the typical and switched configurations in the cases of veno--arterial and veno--venous ECMO with therapeutic plasma exchange. 

In these models, $Q$ is the volume of blood per unit time flowing through the peripheral compartment.  The parameter $\alpha$ is defined so the product $\alpha\, Q$ is the flow through the ECMO compartment and by the TPE device.  The parameter $Q_1$ is the flow into and out of the TPE machine, which is typically much smaller than $Q$.  In all cases, we assume the following constraints on $\alpha$:
\begin{align}
Q_1 < \alpha\, Q < Q \label{eq:alpha1}.
\end{align}  
See Figure \ref{fig:VAtypical}.  These inequalities ensure that the flow into the TPE device does not exceed the flow going by it, and that the ECMO device contributes a fraction $ \alpha$ of the total cardiac output $Q$ through the peripheral compartment.  Define the parameter $[P]_0$ to be the volumetric concentration of total plasma (both new and old plasma) in blood, i.e.~volume of plasma per volume of blood.  We assume $[P]_0$ to be constant during the plasma exchange process, and although it does not explicitly appear in the model equations, we use it for book keeping purposes in the derivation.  The parameters $s_1$, $s_2$, and $s_3$ are transit times through the heart/lung, peripheral, and ECMO compartments, respectively.  In practice, there are distributions of possible transit times, so these parameters can be interpreted as averages of these distributions.  An interesting potential extension for these models would be to incorporate in some way the transit time distributions.  The compartment blood volumes are $V_1$, $V_2$, and $V_3$.  Volumes and transit times are related by the blood flow through each compartment, and these equations will be different in the cases of veno--arterial and veno--venous ECMO.

Referring to Figure \ref{fig:VAtypical}, the variables of interest are $\gamma_1(t)$ and $\gamma_2(t)$, the fractions of new plasma in blood downstream from the heart/lung and peripheral compartments, respectively.  Other auxiliary variables for fraction of new plasma in different locations of the models are used in the derivations and will be defined in context.  We make two important simplifying assumptions that are used in the model derivation for each case. First, we assume old plasma is completely exchanged for new plasma, implying volume of {\em old plasma} per unit time flowing into the device and volume of {\em new plasma} per unit time flowing out of the device is $Q_1[P]_0$. Second, we assume old and new plasma are instantaneously mixed at junctions.  Under these assumptions, conservation of new plasma at junction $A$ in Figures \ref{fig:VAtypical} and \ref{fig:VVtypical}  is: 
\begin{align}
\text{junction }A: \quad (\alpha Q - Q_1)\,\gamma_2(t-s_3)\,[P]_0 + Q_1\,[P]_0 = \alpha Q\,\tilde{\gamma}(t)\,[P]_0, \label{eq:junctionA} 
\end{align}
in which $\tilde{\gamma}(t)$ is the fraction of new plasma directly downstream from the TPE device at time $t$.  The first term in equation \eqref{eq:junctionA} also incorporates the transit time $s_3$ through the ECMO compartment as a delay.  From this point forward, the model derivations for veno--arterial and veno--venous ECMO are different and will be described in the following two subsections. 

\subsection{The typical configuration for veno--arterial ECMO}

\begin{figure}[h!]
\begin{center}
\includegraphics[scale=0.8]{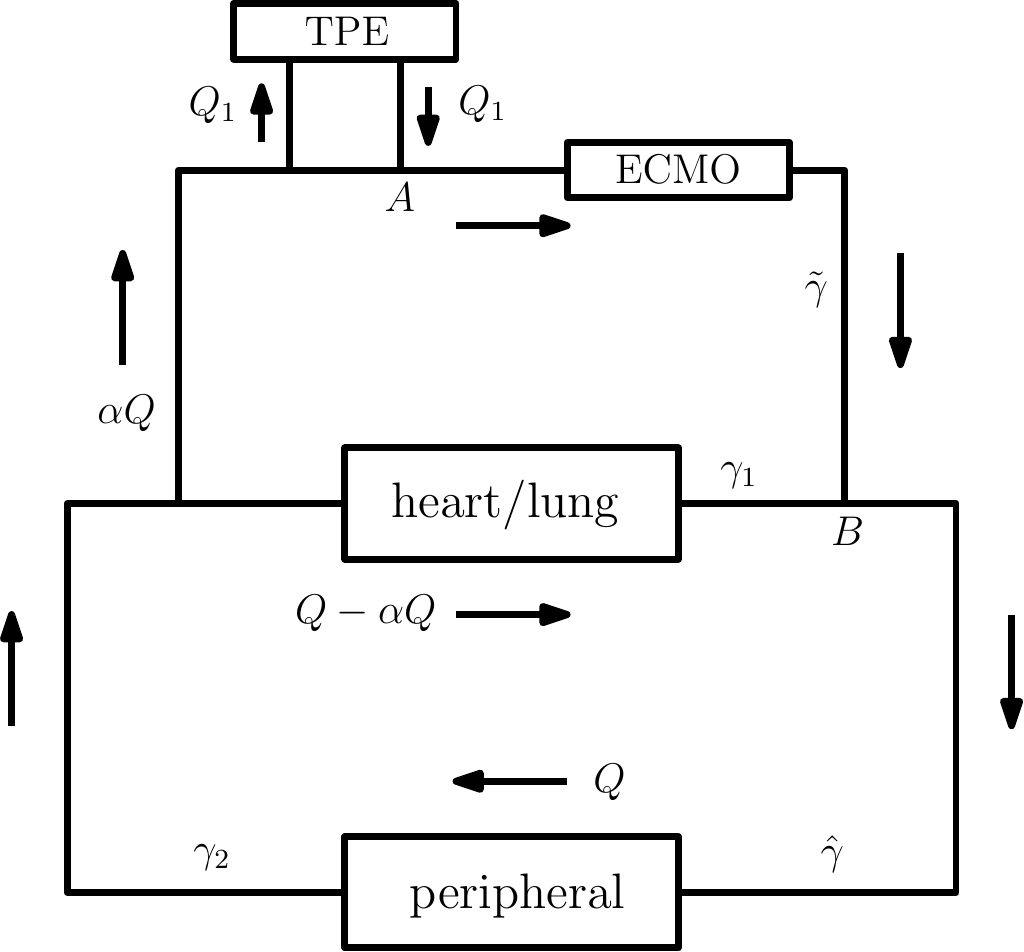}
\caption{The typical configuration for veno--arterial ECMO with therapeutic plasma exchange.}
\label{fig:VAtypical}
\end{center}
\end{figure}

\noindent
For the model of veno--arterial ECMO with TPE, as shown in Figure \ref{fig:VAtypical}, conservation of new plasma at junction $B$ is
\begin{align}
\text{junction }B: \quad \alpha Q\,\tilde{\gamma}(t)\,[P]_0 + (Q - \alpha Q)\,\gamma_1(t)\,[P]_0 = Q\,\hat{\gamma}(t)\,[P]_0, \label{eq:junctionbVAtypical} 
\end{align}
where $\hat{\gamma}(t)$ is the fraction of new plasma directly upstream from the body compartment.  The transit times and volumes of the different compartments are related by the blood flows through each:
\begin{align*}
s_1 = \frac{V_1}{Q - \alpha Q}, \quad s_2 = \frac{V_2}{Q}, \quad s_3 = \frac{V_3}{\alpha Q}.
\end{align*}
Two different models will be derived.  The first model considers the change in volume of new plasma per unit time in the heart/lung and peripheral compartments:
\begin{align*}
\frac{d}{dt}(V_1 \gamma_1) &= (Q - \alpha Q)\, \gamma_2(t) - (Q - \alpha Q)\, \gamma_1(t), \\
\frac{d}{dt}(V_2 \gamma_2) &= Q\, \hat{\gamma}(t) - Q\, \gamma_2(t). 
\end{align*}
Using equations \eqref{eq:junctionA} and \eqref{eq:junctionbVAtypical}, we can rewrite these equations as a system of delay differential equations:
\begin{tcolorbox}
\begin{center}
{\bf DDE model: typical configuration of VA ECMO with TPE}
\end{center}
\begin{equation}
\begin{split}
\frac{d}{dt}(V_1 \gamma_1) &= (Q - \alpha Q)\, \gamma_2(t) - (Q - \alpha Q) \,\gamma_1(t), \\
\frac{d}{dt}(V_2 \gamma_2) &= (Q- \alpha Q)\,\gamma_1(t) + Q_1 + (\alpha Q - Q_1)\,\gamma_2(t-s_3) - Q\, \gamma_2(t),
\end{split}
\label{eq:VAtypicalDDE}
\end{equation}
\end{tcolorbox}

\noindent
which is our first model for the typical configuration of VA ECMO with TPE.  An alternative model takes the form of an algebraic delay equation.  It is derived from equations  \eqref{eq:junctionA}--\eqref{eq:junctionbVAtypical} along with the following relationships between new plasma fractions at different locations and the transit times as delays:
\begin{align*}
\gamma_1(t+s_1) = \gamma_2(t) \quad \text{and} \quad \gamma_2(t+s_2) = \hat{\gamma}(t).
\end{align*}
After algebraic manipulation, the equation for  fraction of new plasma downstream from the heart/lung compartment is:
\begin{tcolorbox}
\begin{center}
\textbf{ADE model: typical configuration of VA ECMO with TPE}
\end{center}
\begin{equation}
\begin{split}
\gamma_1(t) &= \frac{Q_1}{Q}\Big(1  - \gamma_1(t - s_2 - s_3)\Big)\\
&\quad + \alpha\,\gamma_1(t - s_2 - s_3) + (1-\alpha)\,\gamma_1(t - s_1 - s_2).
\label{eq:VAtypicalADE}
\end{split}
\end{equation}
\end{tcolorbox}

\subsection{The typical configuration for veno--venous ECMO}

\begin{figure}[h!]
\begin{center}
\includegraphics[scale=0.8]{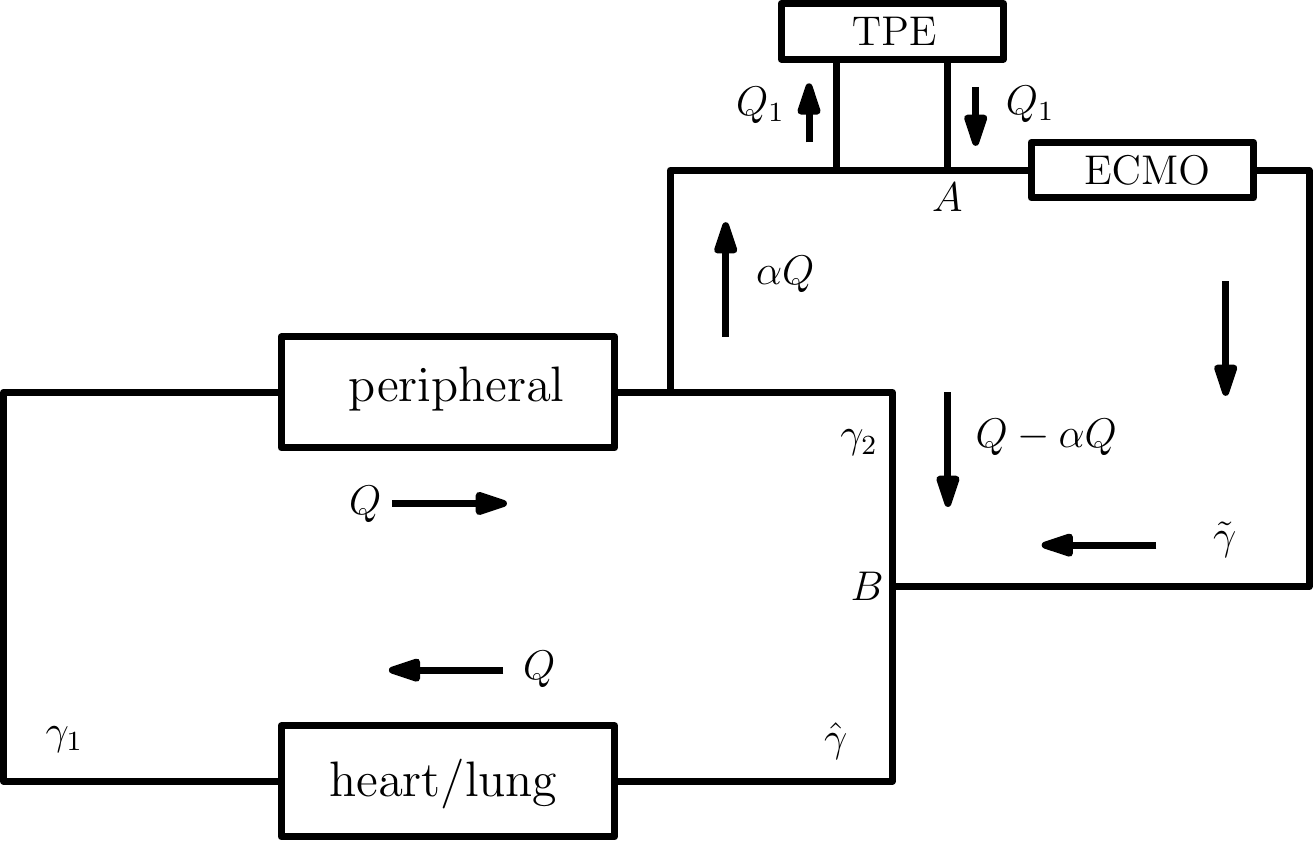}
\caption{Typical configuration for veno--venous ECMO with therapeutic plasma exchange.}
\label{fig:VVtypical}
\end{center}
\end{figure}

\noindent
Figure \ref{fig:VVtypical} depicts our model for veno--venous ECMO with TPE.  As in the previous section, we write down a statement for conservation of new plasma at junction $B$ in Figure \ref{fig:VVtypical}.
\begin{align}
\text{junction }B: \quad \alpha Q\,\tilde{\gamma}(t)\,[P]_0 + (Q-\alpha Q)\,\gamma_2(t)\,[P]_0 = Q\,\hat{\gamma}(t)\,[P]_0. \label{eq:junctionbVVtypical} 
\end{align}
In this case, $\hat{\gamma}(t)$ is the fraction of new plasma upstream from the heart/lung compartment.  As before, the transit times through the different compartments are deteremined by their volumes and corresponding flows:
\begin{align*}
s_1 = \frac{V_1}{Q}, \quad s_2 = \frac{V_2}{Q}, \quad s_3 = \frac{V_3}{\alpha Q}.
\end{align*}
Conservation of new plasma through the peripheral and heart/lung compartments can be written in differential equation form as:
\begin{align*}
\frac{d}{dt}(V_1 \gamma_1) &= Q\, \hat{\gamma}(t) - Q\, \gamma_1(t), \\
\frac{d}{dt}(V_2 \gamma_2) &= Q\, \gamma_1(t) - Q\, \gamma_2(t). 
\end{align*}
Using equations \eqref{eq:junctionA} and \eqref{eq:junctionbVVtypical} we obtain a delay differential equation model:
\begin{tcolorbox}
\begin{center}
{\bf DDE model: typical configuration of VV ECMO with TPE}
\end{center}
\begin{equation}
\begin{split}
\frac{d}{dt}(V_1 \gamma_1) &= (Q-\alpha Q)\,\gamma_2(t) + Q_1 + (\alpha Q - Q_1)\,\gamma_2(t-s_3) - Q\, \gamma_1(t), \\
\frac{d}{dt}(V_2 \gamma_2) &= Q\,\gamma_1(t) - Q\,\gamma_2(t).
\end{split}
\label{eq:VVtypicalDDE}
\end{equation}
\end{tcolorbox}
An algebraic delay equation can also be derived for this configuration, using the following relationships between the transit times and new plasma fractions at different locations:
\begin{align*}
\gamma_1(t+s_1) = \hat{\gamma}(t) \quad \text{and} \quad \gamma_2(t+s_2) = \gamma_1(t).
\end{align*}
After algebraic manipulation we obtain:
\begin{tcolorbox}
\begin{center}
{\bf ADE model: typical configuration of VV ECMO with TPE}
\end{center}
\begin{equation}
\begin{split}
\gamma_1(t) &= \frac{Q_1}{Q}\Big(1  - \gamma_1(t - s_1-s_2 - s_3)\Big) \\
&\quad+ \alpha\gamma_1(t - s_1-s_2 - s_3) + (1-\alpha)\gamma_1(t - s_1 - s_2).
\label{eq:VVtypicalADE}
\end{split}
\end{equation}
\end{tcolorbox}
In this case, the nonzero delay time $s_3$ for the ECMO compartment makes the models above much more interesting.  If transport through the ECMO compartment occurs infinitely fast, i.e. $s_3 = 0$, both models \eqref{eq:VVtypicalDDE} and \eqref{eq:VVtypicalADE} lose their dependence on $\alpha$. 

\subsection{Models for port switching within the plasma exchange device}

\begin{figure}[h!]
\begin{center}
\includegraphics[scale=0.8]{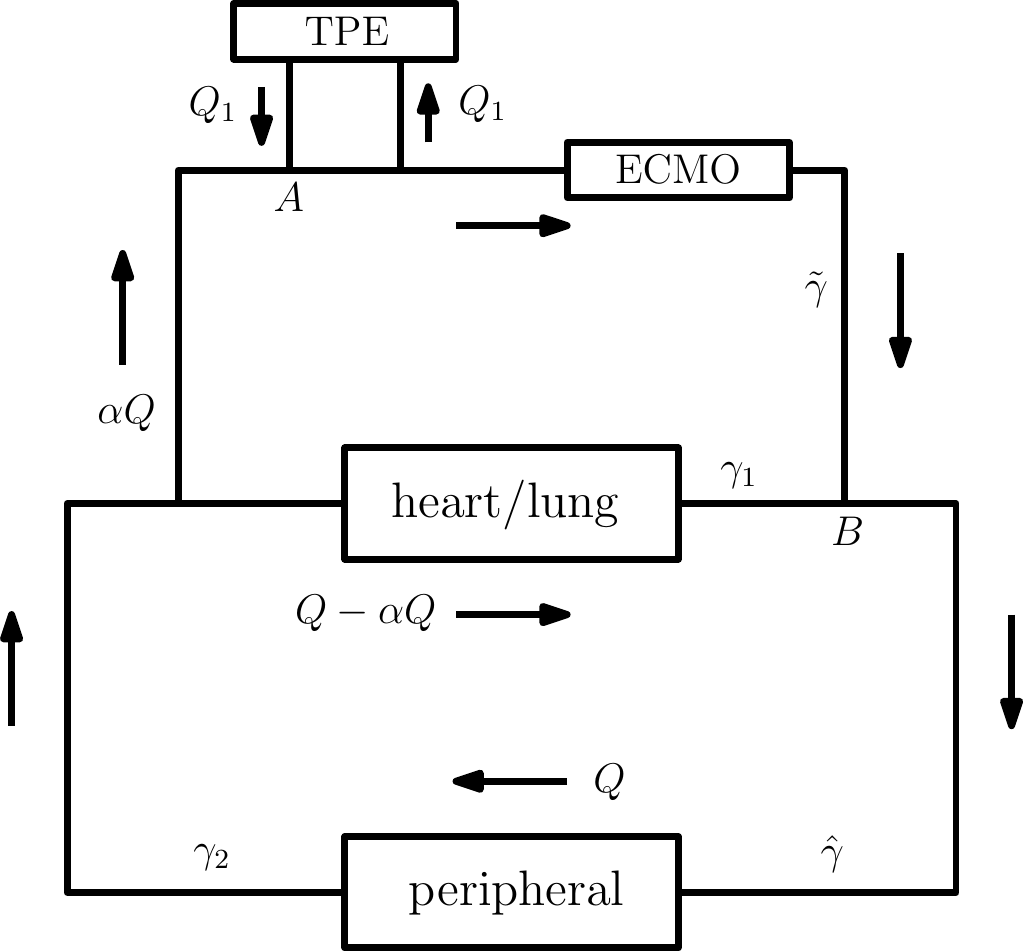}
\caption{The switched configuration for veno--arterial ECMO with therapeutic plasma exchange.}
\label{fig:VAswitched}
\end{center}
\end{figure}

\begin{figure}[h!]
\begin{center}
\includegraphics[scale=0.8]{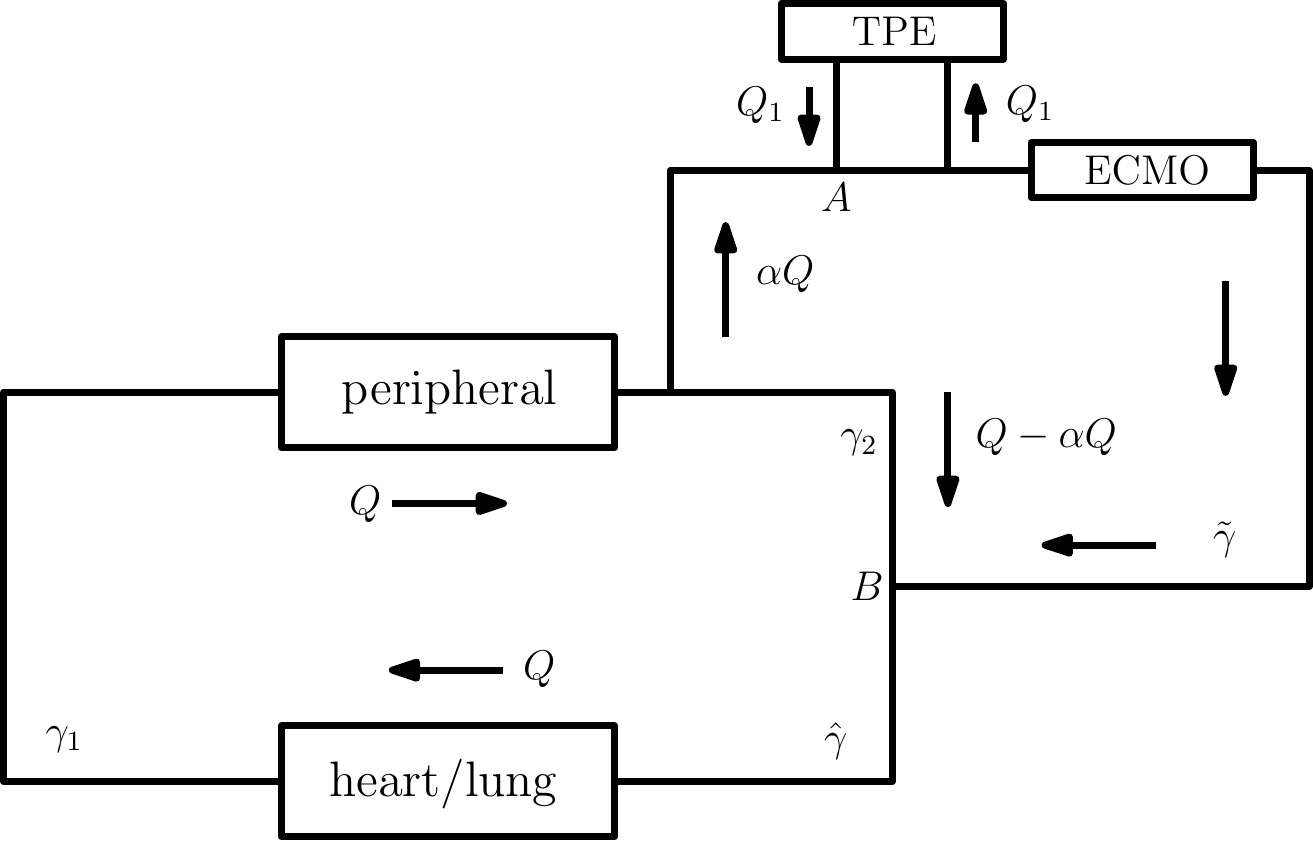}
\caption{The switched configuration for veno--venous ECMO with therapeutic plasma exchange.}
\label{fig:VVswitched}
\end{center}
\end{figure}

\noindent
An application of this modeling approach will be the study of TPE device port switching on the efficiency of the plasma exchange process.  Figures \ref{fig:VAswitched} and \ref{fig:VVswitched} depict schematics for VA and VV ECMO in which the inflow and outflow ports of the TPE device have been reversed from their typical configuration.  In this case, the inflow port is downstream from the outflow port, causing some recirculation of new plasma through the TPE device.

The models for this scenario can be derived in the same way as for the typical configurations.  For veno--arterial ECMO with TPE, we have the following delay differential equation model:
\begin{tcolorbox}
\begin{center}
{\bf DDE model: switched configuration of VA ECMO with TPE}
\end{center}
\begin{equation}
\begin{split}
\frac{d}{dt}(V_1 \gamma_1) &= (Q-\alpha Q)\, \gamma_2(t) - (Q- \alpha Q) \,\gamma_1(t), \\
\frac{d}{dt}(V_2 \gamma_2) &= (Q-\alpha Q)\,\gamma_1(t) + \frac{\alpha Q\, Q_1}{Q_1 + \alpha Q} + \frac{\alpha^2 Q^2}{Q_1 + \alpha Q}\gamma_2(t-s_3) - Q\, \gamma_2(t).
\end{split}
\label{eq:VAswitchedDDE}
\end{equation}
\end{tcolorbox}
\noindent
The algebraic delay equation model in this case is
\begin{tcolorbox}
\begin{center}
{\bf ADE model: switched configuration of VA ECMO with TPE}
\end{center}
\begin{equation}
\begin{split}
\gamma_1(t) &= \frac{\alpha Q_1}{Q_1 + \alpha Q}\Big(1  - \gamma_1(t - s_2 - s_3)\Big)\\
&\quad + \alpha \,\gamma_1(t - s_2 - s_3) + (1-\alpha)\,\gamma_1(t - s_1 - s_2).
\label{eq:VAswitchedADE}
\end{split}
\end{equation}
\end{tcolorbox}
\noindent
For veno--venous ECMO with TPE and the switched configuration for the ports, we have the following delay differential equation model:
\begin{tcolorbox}
\begin{center}
{\bf DDE model: switched configuration of VV ECMO with TPE}
\end{center}
\begin{equation}
\begin{split}
\frac{d}{dt}(V_1 \gamma_1) &= (Q-\alpha Q)\,\gamma_2(t) + \frac{\alpha Q\, Q_1}{Q_1 + \alpha Q} + \frac{\alpha^2 Q^2}{Q_1 + \alpha Q}\gamma_2(t-s_3) - Q\, \gamma_1(t), \\
\frac{d}{dt}(V_2 \gamma_2) &= Q\,\gamma_1(t) - Q\,\gamma_2(t).
\end{split}
\label{eq:VVswitchedDDE}
\end{equation}
\end{tcolorbox}
\noindent
The algebraic delay equation model for this case is:
\begin{tcolorbox}
\begin{center}
{\bf ADE model: switched configuration of VV ECMO with TPE}
\end{center}
\begin{equation}
\begin{split}
\gamma_1(t) &= \frac{\alpha Q_1}{Q_1 + \alpha Q}\Big(1  - \gamma_1(t - s_1-s_2 - s_3)\Big) \\
&\quad+ \alpha\gamma_1(t - s_1-s_2 - s_3) + (1-\alpha)\gamma_1(t - s_1 - s_2).
\label{eq:VVswitchedADE}
\end{split}
\end{equation}
\end{tcolorbox}

\begin{myrem}
Consider a special case for the models of VV ECMO with TPE in which the heart/lung and peripheral compartments are treated as a single compartment.  In this case, assume the fractions of new plasma downstream from each of these compartments are the same:
\begin{align*}
\gamma_1(t) = \gamma_2(t).
\end{align*} 
 Also, assume the transit time through the ECMO compartment is zero, corresponding to $s_3 = 0$.  These assumptions formally convert our model to the single compartment one for plasma exchange of Randerson et al. \cite{Randerson82}. Note that the DDE and ADE models for the typical configuration of VV ECMO with TPE do not depend on $\alpha$.  Interestingly, for the switched configuration, the models still depend on $\alpha$.  
The DDE models for the typical and switched configurations for VV ECMO with TPE take the simplified form:
\begin{align}
\frac{d}{dt} \gamma_2 = \frac{\beta}{V_1 + V_2} \left( 1 - \gamma_2(t) \right),
\label{eq:analyticVV}
\end{align}
where for the typical configuration we have $\beta = Q_1$ and for the switched configuration we have $\beta = \frac{ \alpha Q\, Q_1 }{ \alpha Q + Q_1}.$ Equation \eqref{eq:analyticVV} for $\gamma_2$ has an analytical solution of the form:
\begin{align*}
\gamma_2(t) = 1 - \exp\left(-\frac{\beta\, t}{V_1 + V_2} \right).
\end{align*}
For the typical configuration, this equation can be rewritten as:
\begin{align}
1-\gamma_2(t) = \exp\left(-\frac{Q_1[P]_0\,t}{(V_1 + V_2)[P]_0} \right).
\label{eq:analyticVV2}
\end{align}
The term $1 - \gamma_2(t)$ is the fraction of old plasma downstream from the heart/lung and peripheral compartments, at time $t$.  The term $Q_1[P]_0\, t$ is the volume of plasma processed up to time $t$, since $Q_1$ is the TPE device flow, and $(V_1 + V_2)[P]_0$ is the total volume of plasma, both new and old, within the native circulation (i.e. the sum of volumes in the heart/lung and peripheral compartments).  The ratio of these two terms appearing in the exponential of \eqref{eq:analyticVV2} can be interpreted as the {\em number of plasma volumes processed}, corresponding to a multiple of the total volume of plasma contained in the native circulation.  In words, equation \eqref{eq:analyticVV2} is:
\begin{tcolorbox}
\begin{center}
{\bf typical configuration}
\end{center}
\begin{align}
\left(
\begin{gathered}
\text{fraction of old} 
\\ \text{plasma remaining}
\end{gathered} \right) = \exp\left( -\text{plasma volumes processed} \right),
\label{eq:analyticVV3}
\end{align} 
\end{tcolorbox}
\noindent
which is one of the ``informal laws of apheresis.''  Equation \eqref{eq:analyticVV3} is the same as that obtained in the single compartment model derived by Randerson et al. \cite{Randerson82}, with the seiving coefficient of the plasma exchange device equal to 1.  This is significant since their model of plasma exchange was validated with clinical data.  The analogous equation for the switched configuration is: 
\begin{tcolorbox}
\begin{center}
{\bf switched configuration}
\end{center}
\begin{align}
\left(
\begin{gathered}
\text{fraction of old} 
\\ \text{plasma remaining}
\end{gathered} \right) = \exp\left( - \frac{\alpha Q}{\alpha Q + Q_1} \times \text{plasma volumes processed} \right).
\label{eq:analyticVV4}
\end{align} 
\end{tcolorbox}
\noindent
Equation \eqref{eq:analyticVV4} is slight modification of \eqref{eq:analyticVV3} and accounts for port switching in the plasma exchange device; it includes the term $\frac{\alpha Q}{\alpha Q + Q_1}$ that multiplies the number of plasma volumes.  Practically, the TPE device flow $Q_1$ is small relative to $\alpha Q$, so the ratio $\frac{\alpha Q}{\alpha Q + Q_1} \approx 1$ and the typical and switched configurations have roughly the same plasma exchange efficiency.  As we shall show in our numerical results, the switched configuration is much less efficient in the regime $\alpha Q \approx Q_1$.
\end{myrem} 

\begin{myrem}
\label{rem:rem2}
The algebraic delay equations can be derived from the delay differential equations by replacing the time derivative with an approximate difference quotient.  In particular, the difference quotient
\begin{align}
\frac{d}{dt}(V_i \gamma_i) \approx \frac{V_i}{s_i}\left(\gamma_i(t + s_i) - \gamma_i(t)\right), \quad i = 1,2,
\label{eq:dq}
\end{align}
used in the DDE model, combined with the relationships between the flows, volumes, and transit times, results in the corresponding ADE model.  As an example, equation \eqref{eq:VAtypicalDDE} with approximation \eqref{eq:dq} becomes 
\begin{align}
\frac{V_1}{s_1}\left(\gamma_1(t + s_1) - \gamma_1(t)\right)  &= (Q - \alpha Q)\, \gamma_2(t) - (Q - \alpha Q) \,\gamma_1(t), \label{eq:simp1} \\
\frac{V_2}{s_2}\left(\gamma_2(t + s_2) - \gamma_2(t)\right) &= (Q- \alpha Q)\,\gamma_1(t) + Q_1 + (\alpha Q - Q_1)\,\gamma_2(t-s_3) - Q\, \gamma_2(t). \label{eq:simp2}
\end{align}
Using $(Q - \alpha Q)\,s_1 = V_1$ in equation \eqref{eq:simp1}, we obtain the delay relation used to derive the algebraic delay equation, namely:
\begin{align*}
\gamma_1(t+s_1) = \gamma_2(t).
\end{align*}
The algebraic delay equation \eqref{eq:VAtypicalADE} for this case follows from the above equation along with \eqref{eq:simp2} and $Q\,s_2 = V_2$.  Similar arguments can be made for the other configurations.  

In this way, the ADE models can be interpreted as forward Euler discretizations of the DDE models, with timesteps equal to compartment transit times. Depending on the clinical scenario, plasma exchange may be done multiple times, and each exchange procedure happens over the course of hours.  The time scale of this process is slow in comparison to transit times through various compartments of the body. For example, if the volume of blood is 5 liters and the corresponding cardiac output is 5 liters/minute, the transit time of a red blood cell through the entire circulation can be approximated as the ratio of the volume to the flow, i.e. 60 seconds.  This separation of time scales reveals that the ADE and DDE models provide results very close to each other with a difference on the order of $\max(s_1,s_2)$.   
\end{myrem}

\section{Results}
\label{sec:results}

In this section, we describe some results from numerical simulations of our models for veno--arterial and veno--venous ECMO with plasma exchange.  First, we compare the algebraic delay equations and delay differential equations.  Second, we use these models to quantify the differences between the typical and switched configurations of the plasma exchange device.  Third, we perform an analysis of the sensitivity of these models to their parameters. In our simulations, the timestep size is chosen to be $\Delta t = 10^{-2}$ seconds.  The transit times $s_1$, $s_2$, and $s_3$ are chosen in a manner to be described below based on parameters and other available data, but are selected to be the closest multiple of the timestep size for ease in simulating the models.  The delay differential equations are discretized with the forward Euler method.  

\subsection{Comparing algebraic delay and delay differential equation models}
\label{subsec:1}

First, we compare the ADE and DDE models for each configuration.  Nominal parameters are chosen as follows.  The flow $Q$ is set to 116.7 mL/second \cite{Weiss09}.  Transit times through the peripheral and heart/lung compartments are chosen to be $s_1 = 13$ seconds and $s_2 = 39$ seconds, following estimations of compartment volumes in Weiss and Neimann et al. \cite{Weiss09, Niemann02}.  We assume the ECMO circuit maintains 70\% of the flow through the peripheral compartment corresponding to $\alpha = 0.7$.  The transit time through the ECMO compartment is determined by $s_3 = \frac{V_3}{\alpha Q}$, where the ECMO blood volume $V_3$ is chosen to be 500 mL \cite{Belousova19, Jhang07}.  The TPE device flow $Q_1$ is set to 1.5 mL/second \cite{Dyer14}.  Initial conditions for fractions of new plasma are $\gamma_1(t) = \gamma_2(t) = 0$ for $0 \leq t \leq s_1 + s_2 + s_3$.

\begin{figure}[h!]
\begin{center}
\includegraphics[scale=0.45,trim=0 0 -50 0]{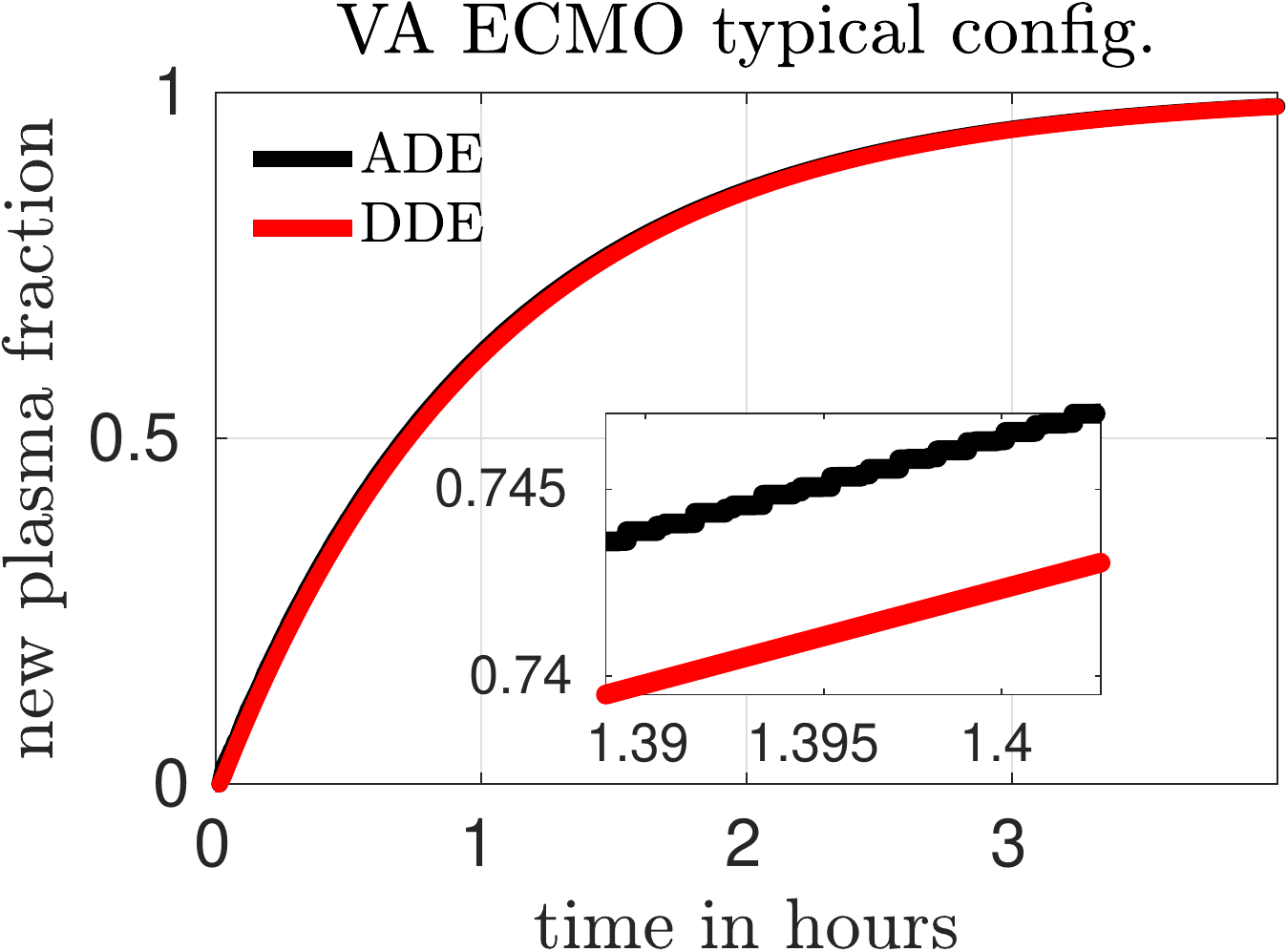}
\includegraphics[scale=0.45]{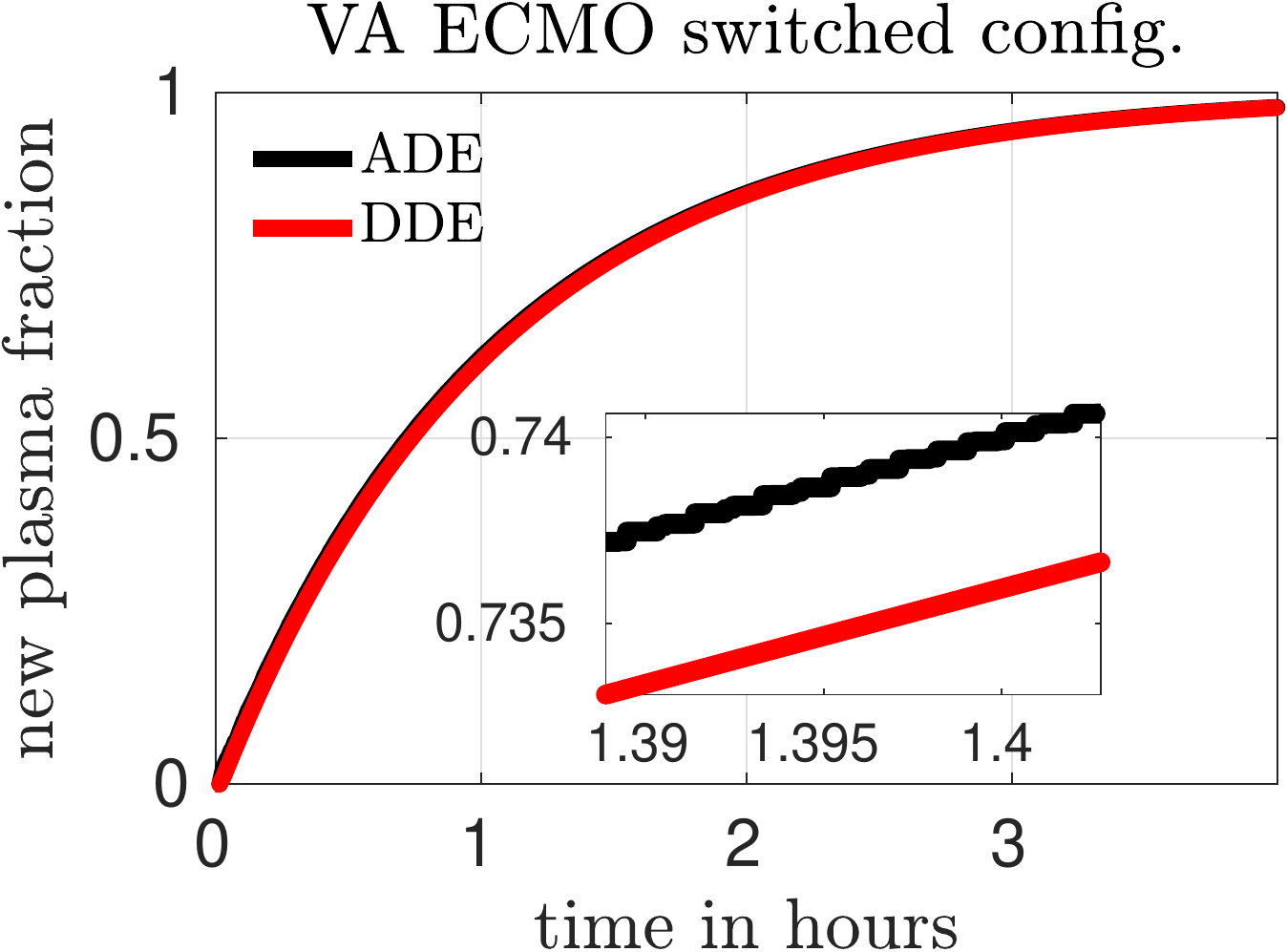}
\caption{A comparison of the ADE (black) and DDE (red) models for VA ECMO with TPE, for the typical configuration of the TPE on the left, and the switched configuration on the right.  The inset figure shows a zoomed in portion of the results.}
\label{fig:compareVA}
\end{center}
\end{figure}

\begin{figure}[h!]
\begin{center}
\includegraphics[scale=0.45,trim=0 0 -50 0]{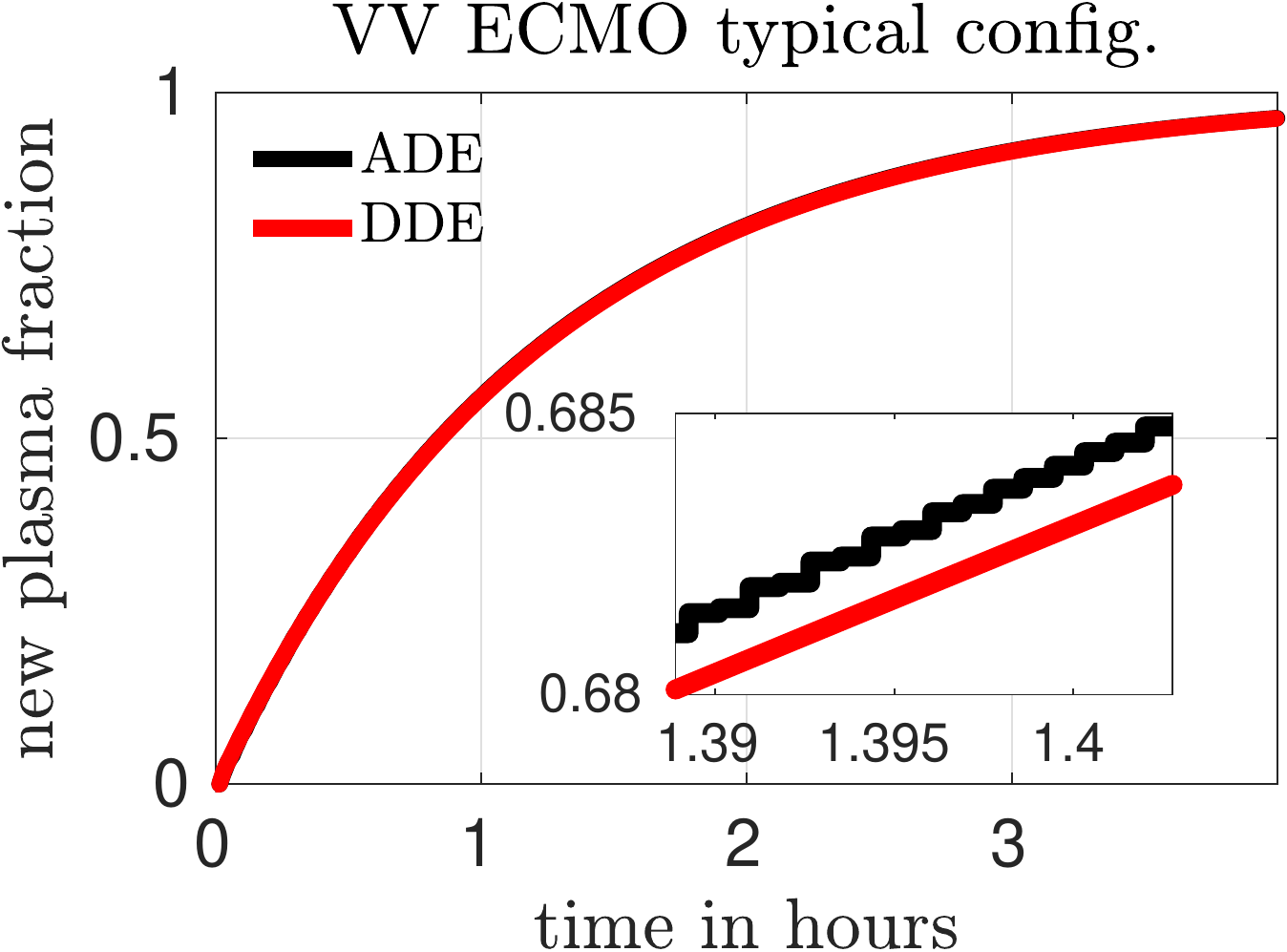}
\includegraphics[scale=0.45]{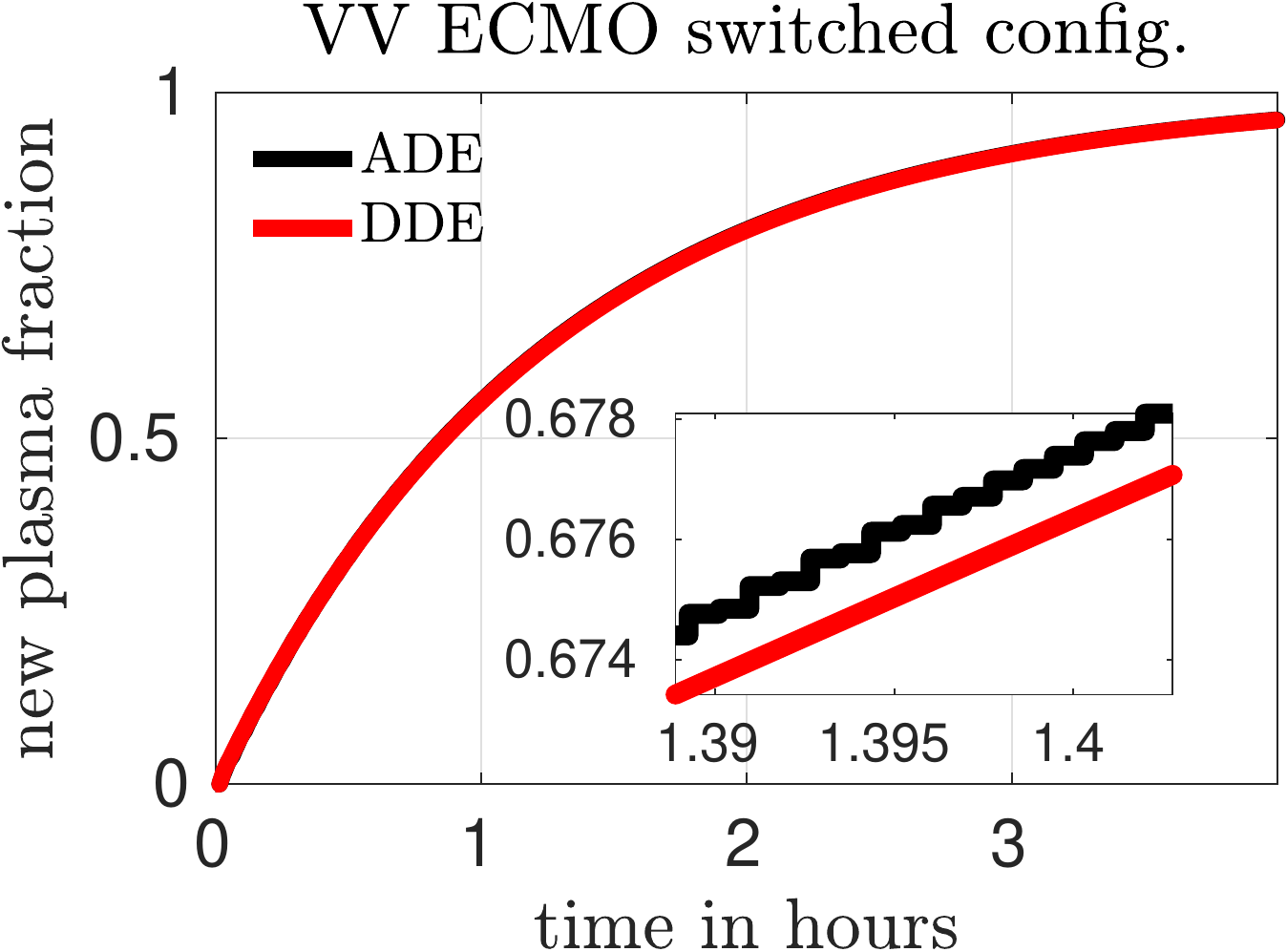}
\caption{A comparison of the ADE (black) and DDE (red) models for VV ECMO with TPE, for the typical configuration of the TPE on the left, and the switched configuration on the right.  The inset figure shows a zoomed in portion of the results.}
\label{fig:compareVV}
\end{center}
\end{figure}

Figure \ref{fig:compareVA} shows a comparison of the ADE and DDE models for VA ECMO.  The new plasma fraction downstream from the heart/lung compartment, $\gamma_1(t)$, is plotted for each case.  Results for the typical configuration of the plasma exchange device are on the left and results for the switched configuration are on the right.  Figure \ref{fig:compareVV} shows the corresponding results for VV ECMO.  In both sets of figures we  included an inset figure which depicts a small subinterval of time to show the small differences between the ADE and DDE models.  Following Remark \ref{rem:rem2}, the two types of models are very close to each other since the compartment transit times are on the order of minutes, and the plasma exchange procedure is simulated over the course of multiple hours.  We also remark that the inset figures show some small differences between the typical and switched configurations.  The switched configuration is less efficient since the fraction of new plasma at a given time is smaller than in the typical configuration.  We further investigate these differences in the next subsection.

\subsection{The effect of port switching on plasma exchange}

In this subsection, we use our models to compare the typical and switched configurations for the TPE device ports.  In the switched configuration, the return line for the TPE device is upstream from the inlet line.  This setup leads to recirculation of new plasma, meaning that some new plasma entering the ECMO circuit through the return line is immediately processed by the TPE device.  The switched configuration is not an ideal setup for plasma exchange, but our calculations in Remark \ref{rem:rem2} and simulations from the previous subsection suggest that it results in very similar plasma exchange efficiency for a nominal choice of parameters.  As also suggested in Remark \ref{rem:rem2}, the switched configuration may be much less efficient when the fraction $\alpha$ of cardiac output supported by the ECMO device is very small.

\begin{figure}[h!]
\begin{center}
\includegraphics[scale=0.45,trim=0 0 -50 0]{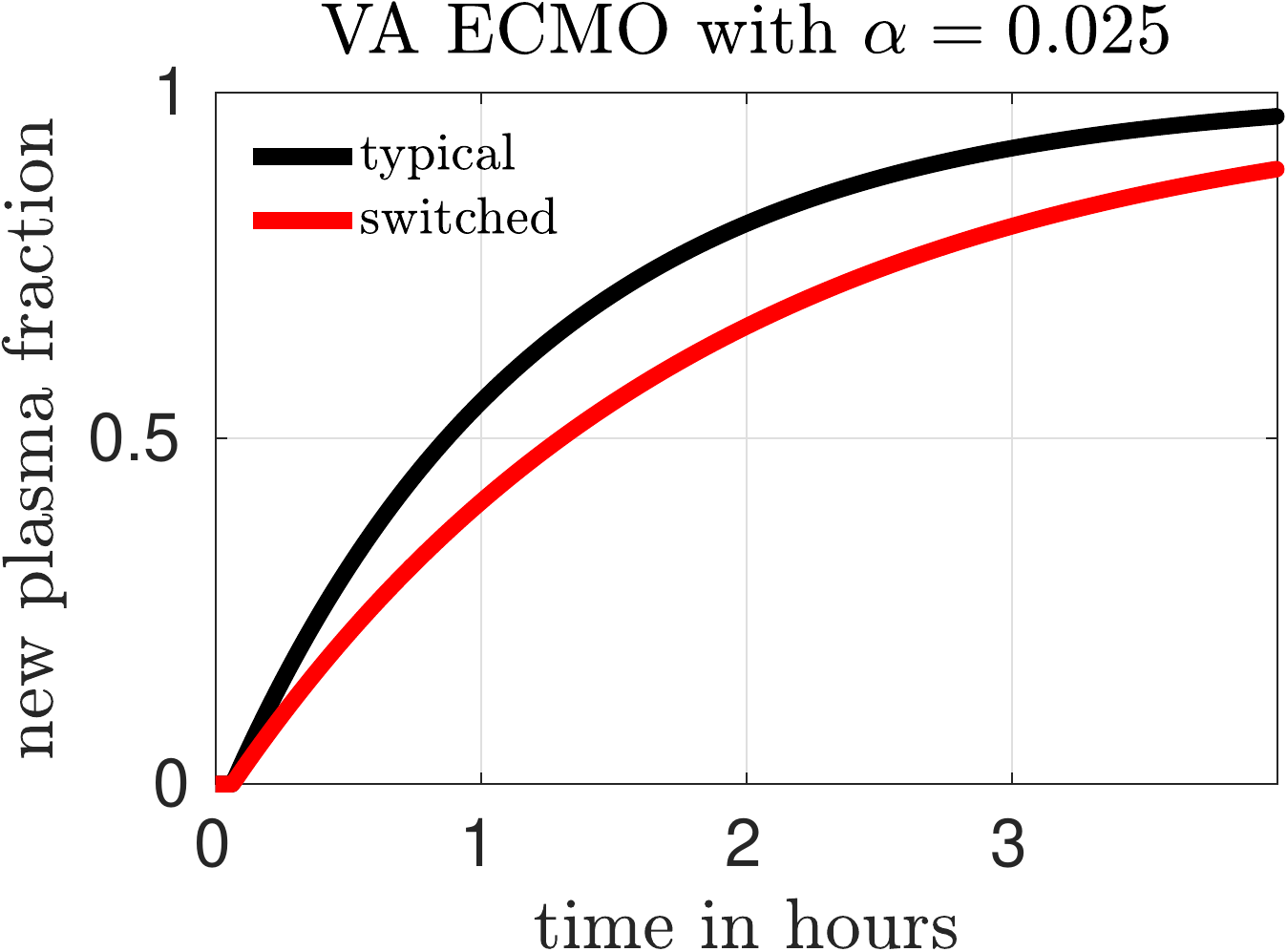}
\includegraphics[scale=0.45]{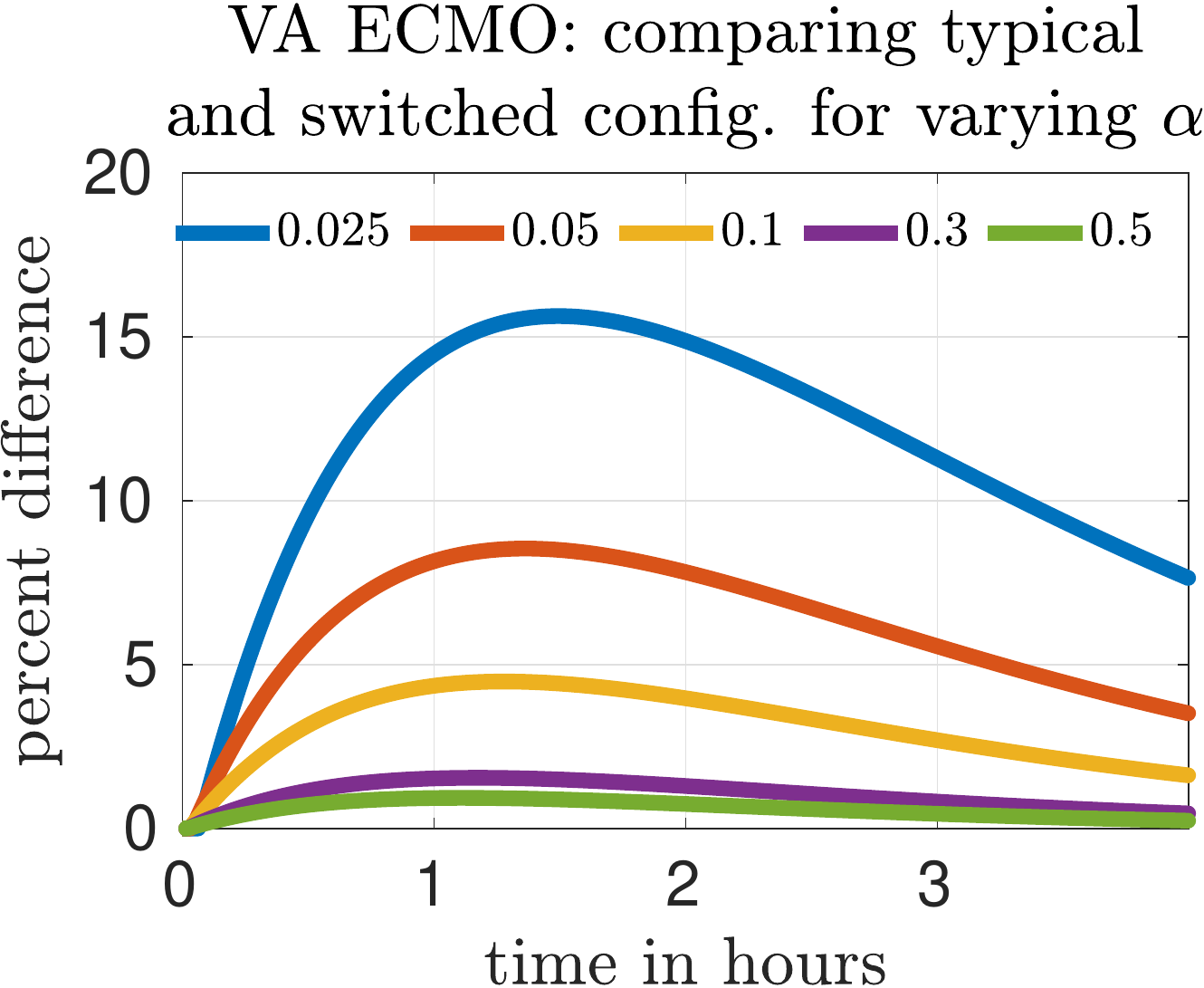}
\caption{Results for VA ECMO: the left panel shows fraction of new plasma for a small value of $\alpha$ for both the typical and switched configurations. The right panel shows the percent difference in fraction of new plasma between these two configurations for several values of $\alpha$.}
\label{fig:switchingVA}
\end{center}
\end{figure}

\begin{figure}[h!]
\begin{center}
\includegraphics[scale=0.45,trim=0 0 -50 0]{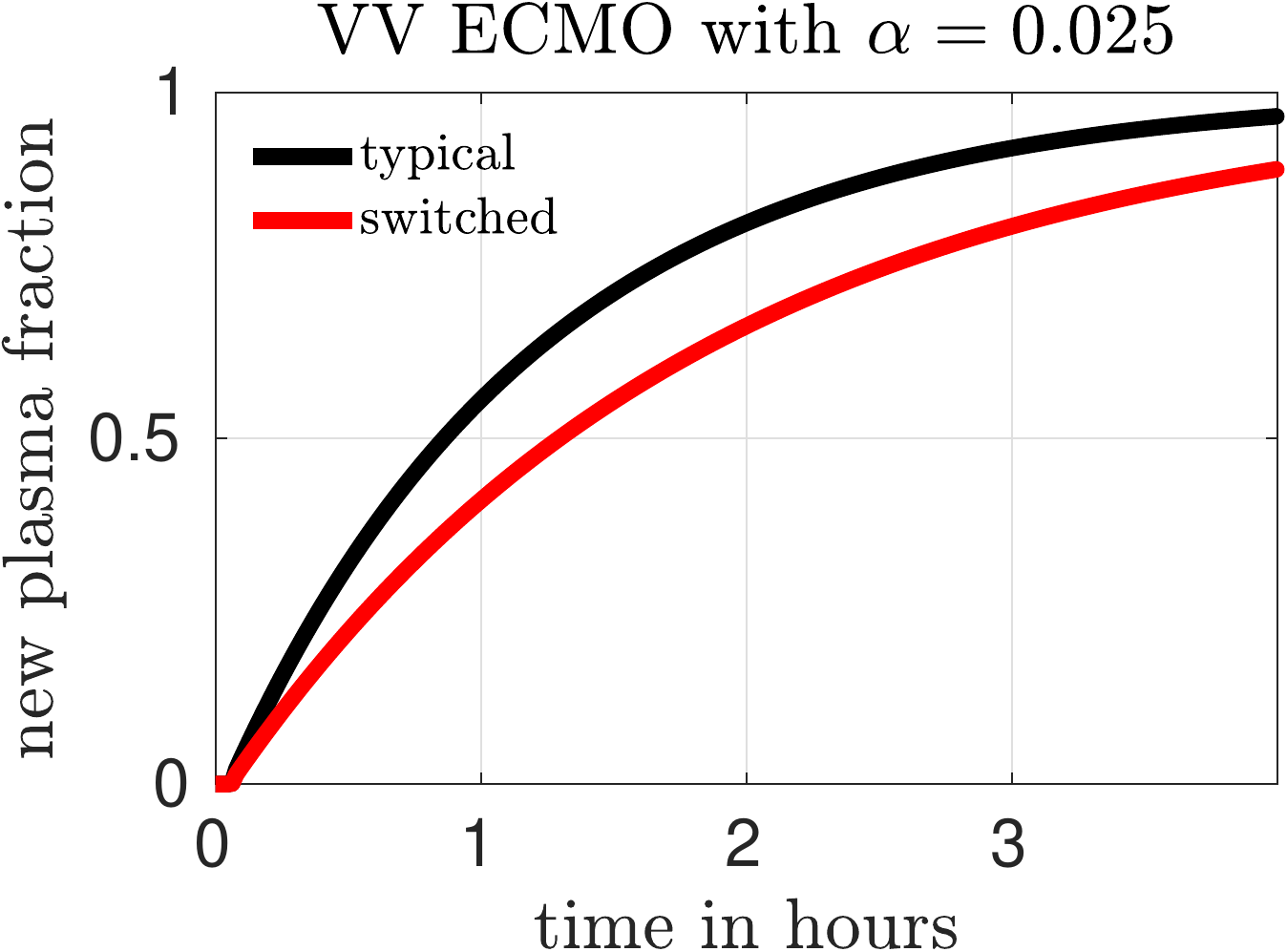}
\includegraphics[scale=0.45]{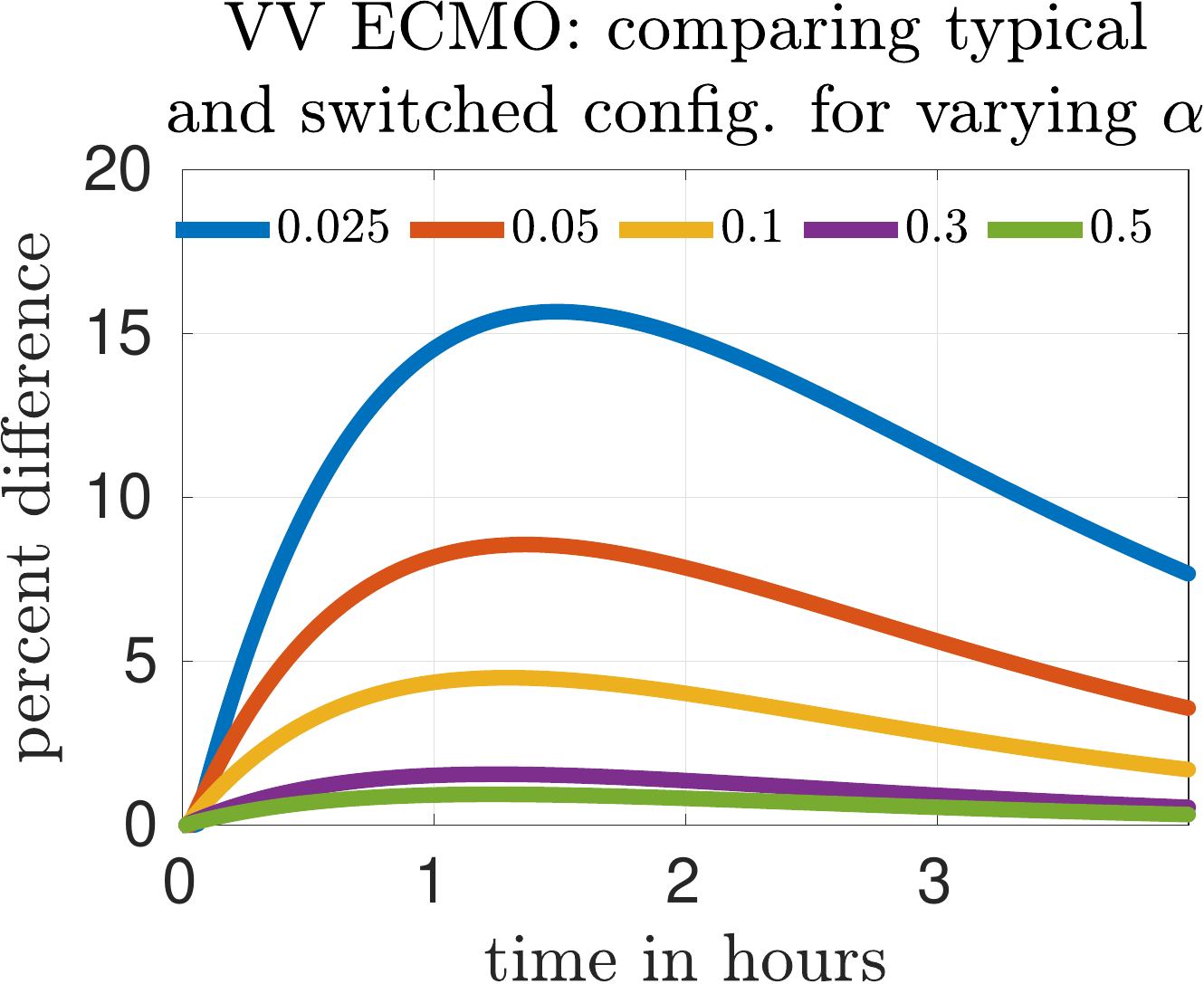}
\caption{Results for VV ECMO: the left panel shows fraction of new plasma for a small value of $\alpha$ for both the typical and switched configurations. The right panel shows the percent difference in fraction of new plasma between these two configurations for several values of $\alpha$.}
\label{fig:switchingVV}
\end{center}
\end{figure}

Figures \ref{fig:switchingVA} and \ref{fig:switchingVV} depict results from our models with smaller values of $\alpha$ and provide a comparison of the typical and switched configurations.  Simulations are done with the DDE models and the same nominal parameters as in the previous section.  In the left panel of these figures, we show new plasma fraction curves for both configurations and for the smallest value of $\alpha$ considered.  The switched configuration for the cases of VA and VV ECMO is much less efficient in exchanging plasma than the typical configuration.  The right panel shows percent differences between these configurations for different values of $\alpha$. Notice that as $\alpha$ grows, both configurations maintain essentially the same plasma exchange efficiency.

\subsection{Sensitivity analysis}

This subsection describes a study of the sensitivity of the fraction of new plasma to changes in the parameters of the model.  We define the sensitivity for a parameter $Y \in \{Q_1, Q, \alpha, s_1, s_2, V_3 \}$ as the partial derivative of $\gamma_1$ with respect to $Y$, evaluated at $t = $ 4 hours.  The sensitivity is approximated as:
\begin{align*}
\frac{\partial \gamma_1}{\partial Y}\Big|_{t = 4 \text{ hours}} \approx \frac{\gamma_1(t,\tilde{Y}) - \gamma_1(t,Y)}{\tilde{Y} - Y}\Big|_{t = 4 \text{ hours}}, 
\end{align*}
in which $Y$ is taken to be the nominal value of the parameter used in subsection \ref{subsec:1}, and $\tilde{Y} = 1.1 \times Y$, i.e. ten percent larger than the nominal value.  In other words, our approximation to the sensitivity is a forward difference approximation to the partial derivative.  For the calculations in this section, the ADE models for both the typical and switched configuration are used. 

\begin{figure}[h!]
\begin{center}
\includegraphics[scale=0.45,trim=0 0 -50 0]{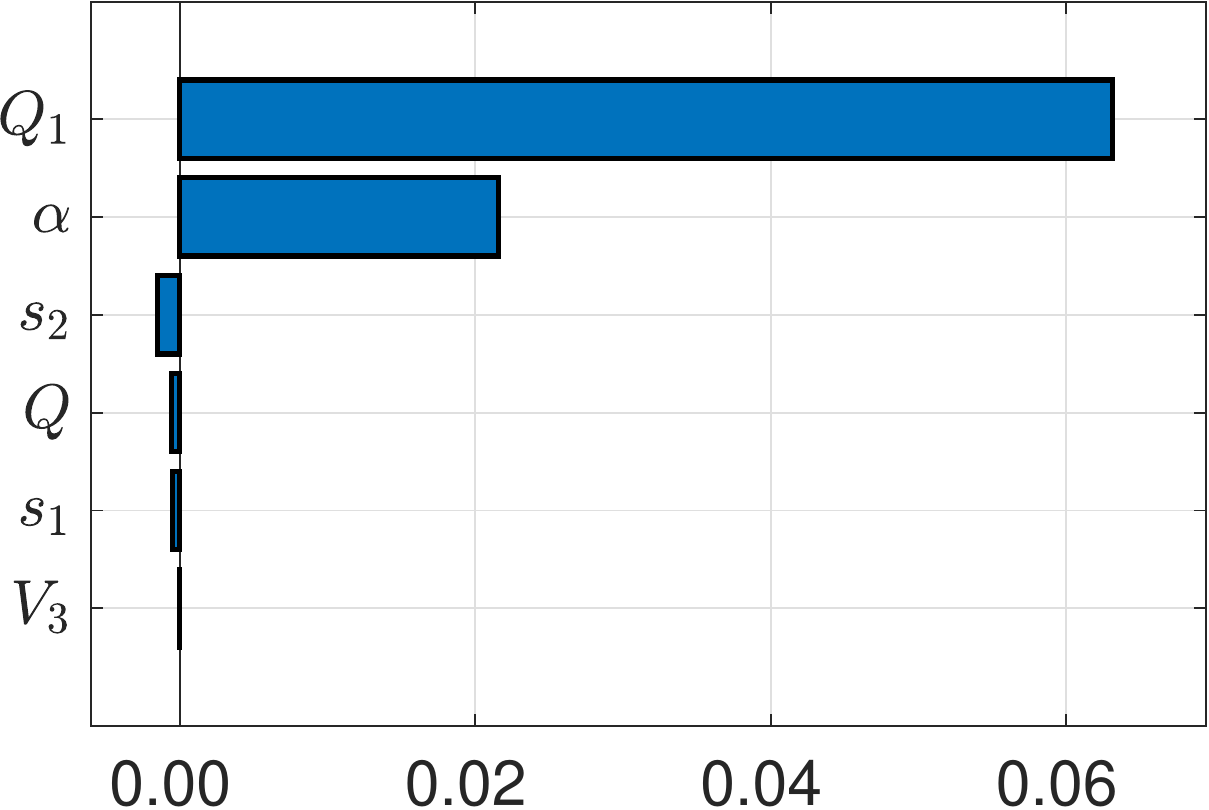}
\includegraphics[scale=0.45]{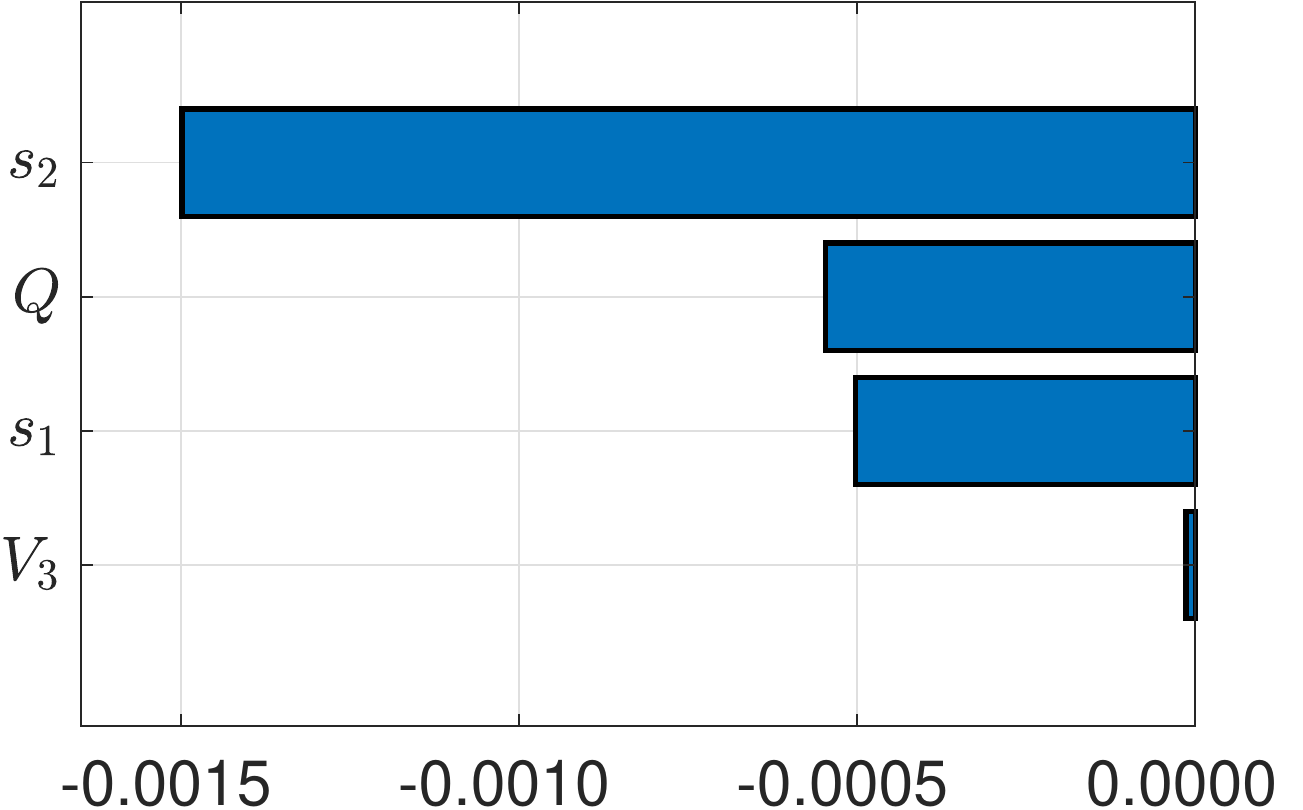}
\put(-310,112){Sensitivities for VA ECMO -- typical config. -- ADE model}
\caption{Sensitivities for VA ECMO with the typical configuration of the TPE device.  The left panel shows all of the sensitivities, and the right panel shows the four smallest sensitivities.}
\label{fig:sensVAtypicalADE}
\end{center}
\end{figure}

\begin{figure}[h!]
\begin{center}
\includegraphics[scale=0.45,trim=0 0 -50 0]{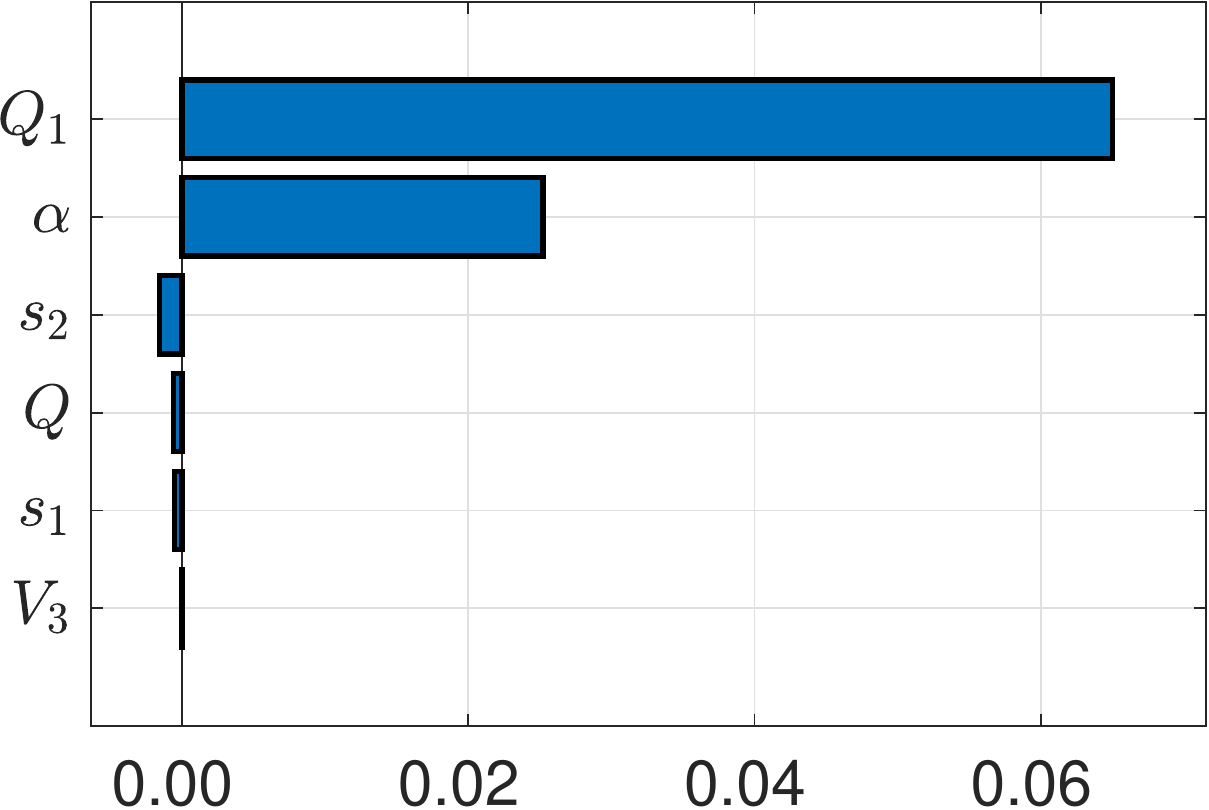}
\includegraphics[scale=0.45]{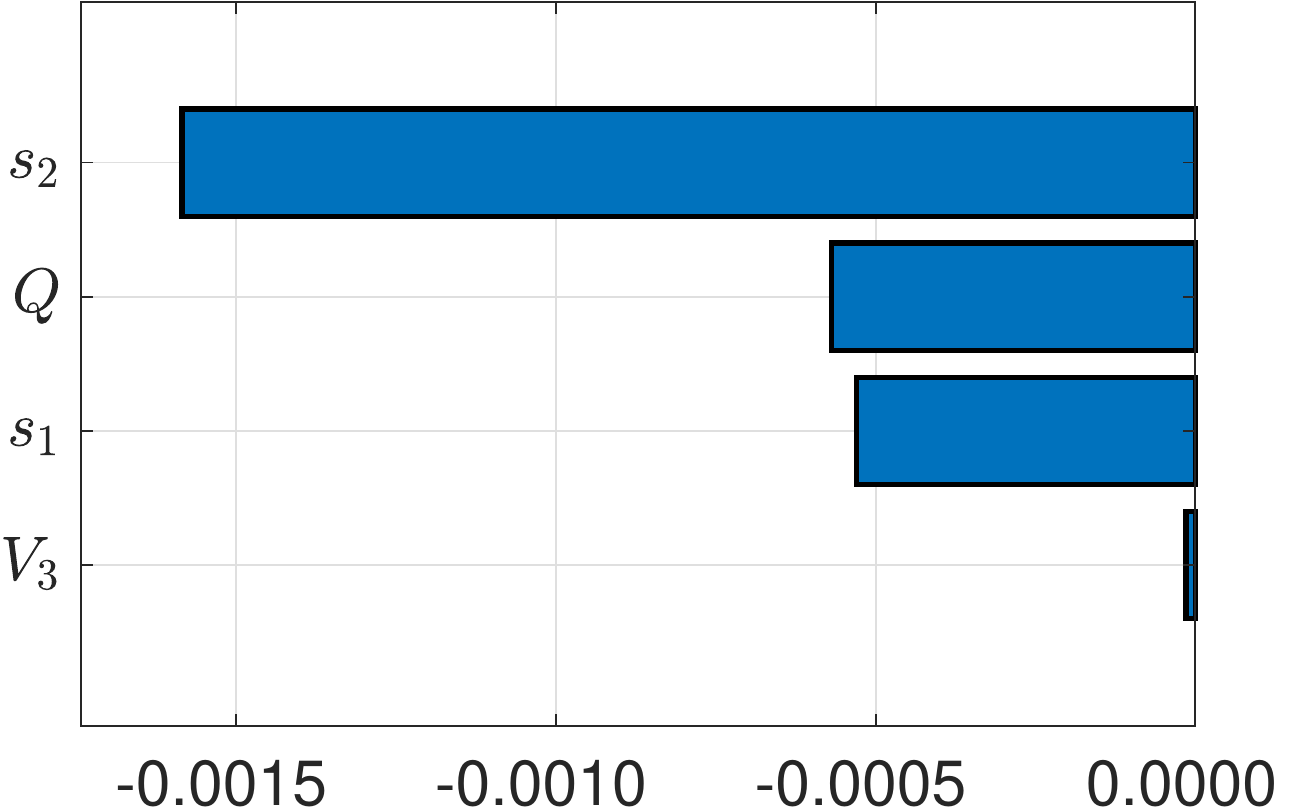}
\put(-308,112){Sensitivities for VA ECMO -- switched config. -- ADE model}
\caption{Sensitivities for VA ECMO with the switched configuration of the TPE device.  The left panel shows all of the sensitivities, and the right panel shows the four smallest sensitivities.}
\label{fig:sensVAswitchedADE}
\end{center}
\end{figure}

\begin{figure}[h!]
\begin{center}
\includegraphics[scale=0.45,trim=0 0 -50 0]{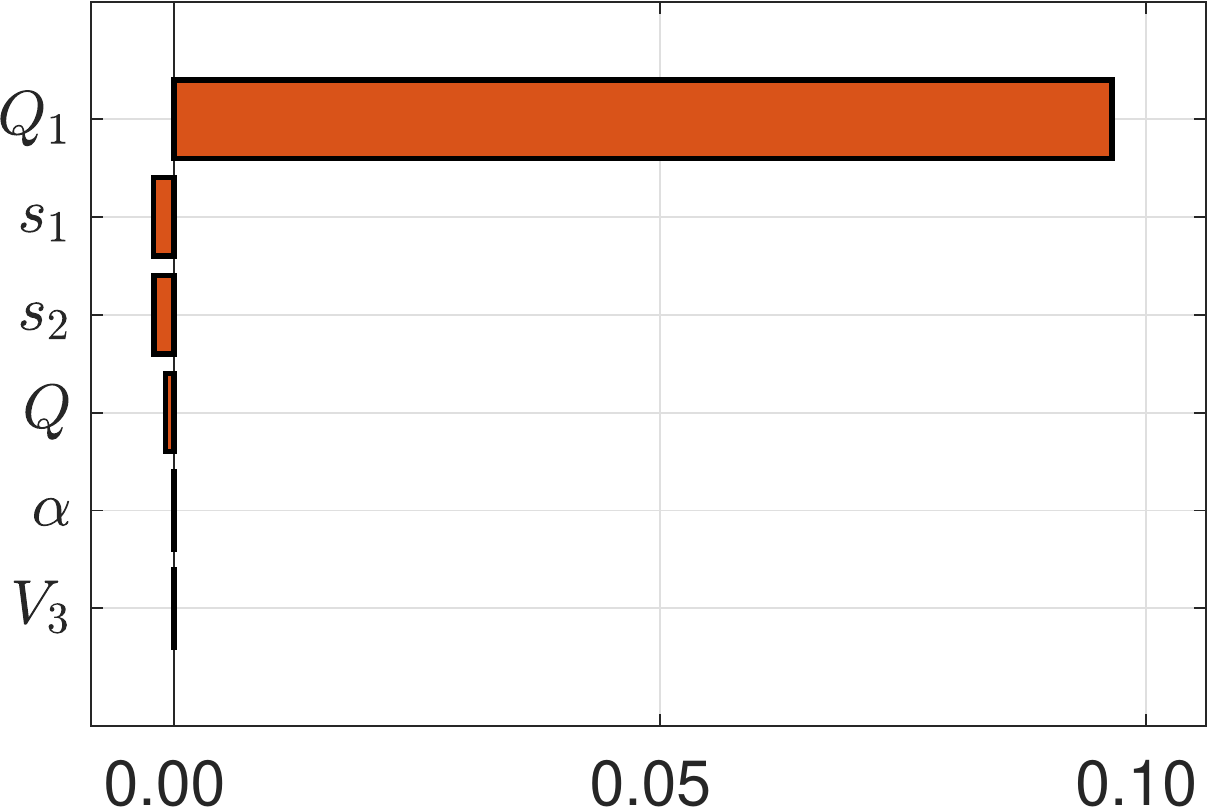}
\includegraphics[scale=0.45]{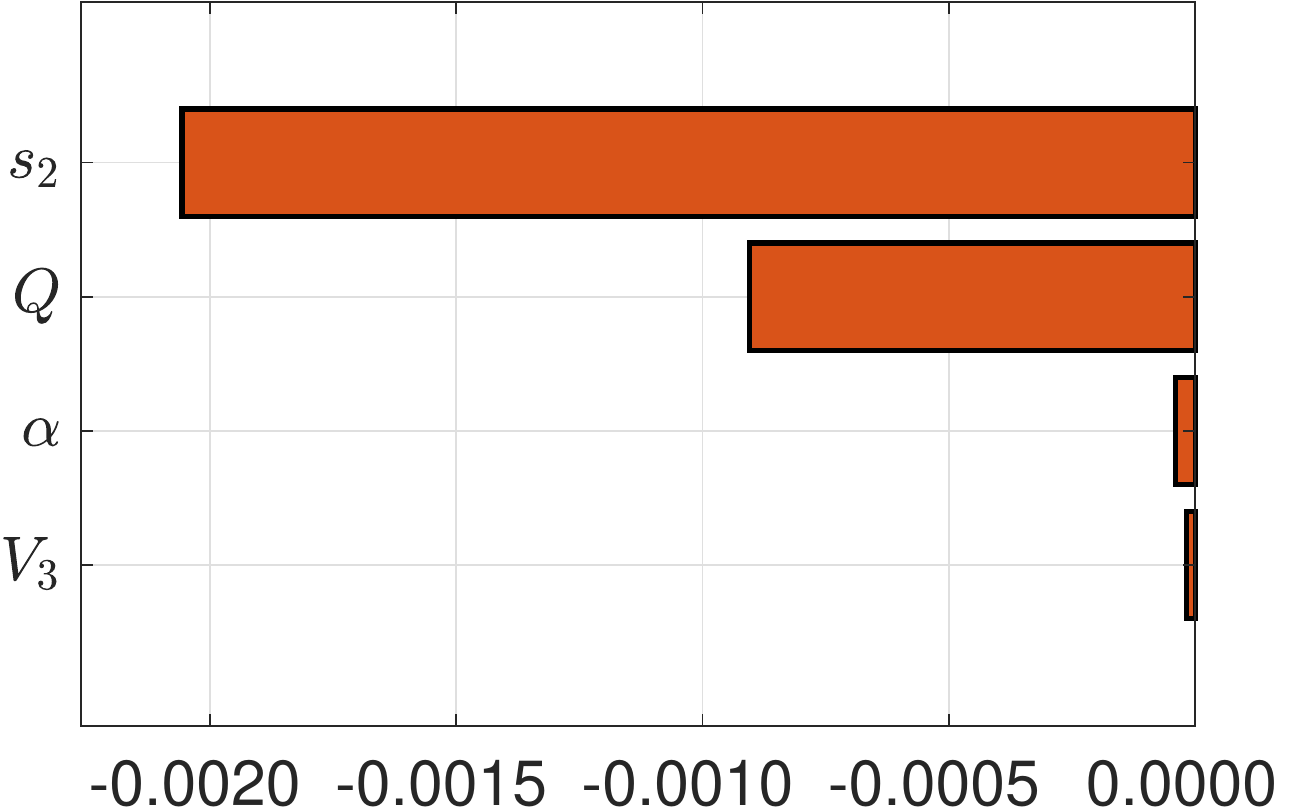}
\put(-310,112){Sensitivities for VV ECMO -- typical config. -- ADE model}
\caption{Sensitivities for VV ECMO with the typical configuration of the TPE device.  The left panel shows all of the sensitivities, and the right panel shows the four smallest sensitivities.}
\label{fig:sensVVtypicalADE}
\end{center}
\end{figure}

\begin{figure}[h!]
\begin{center}
\includegraphics[scale=0.45,trim=0 0 -50 0]{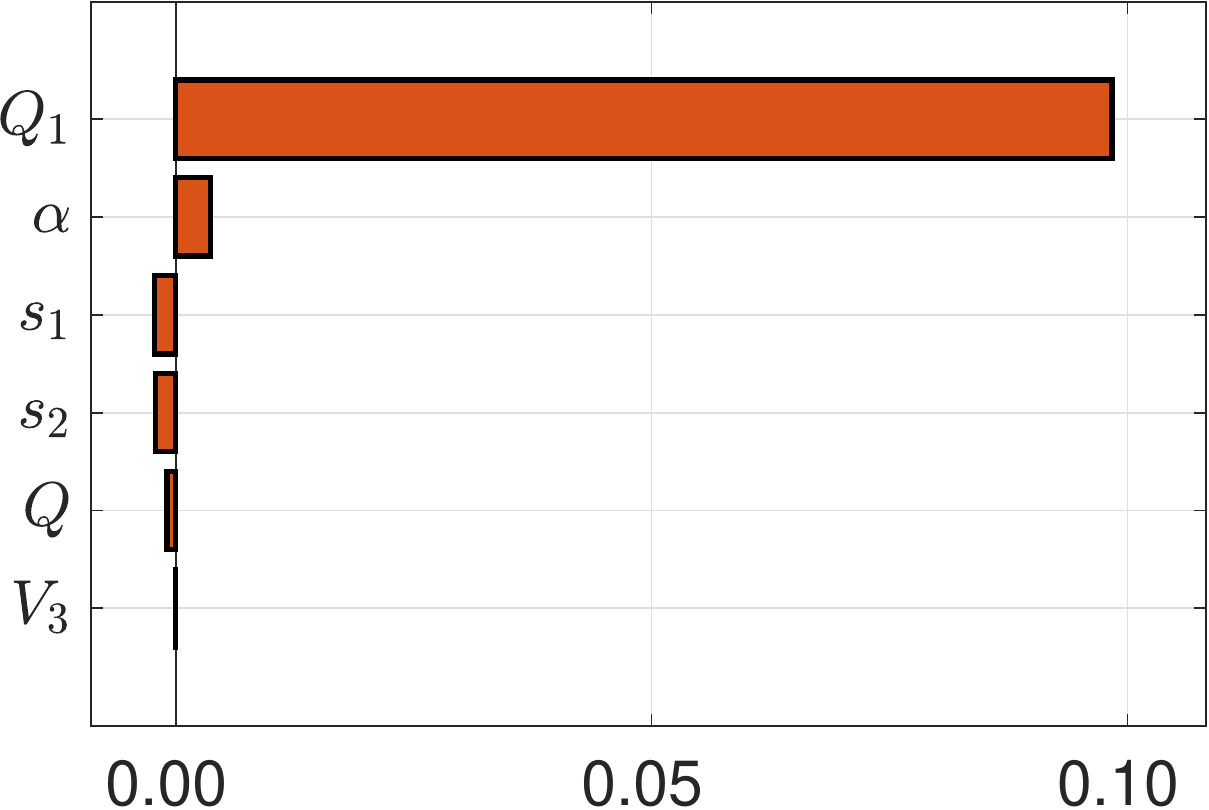}
\includegraphics[scale=0.45]{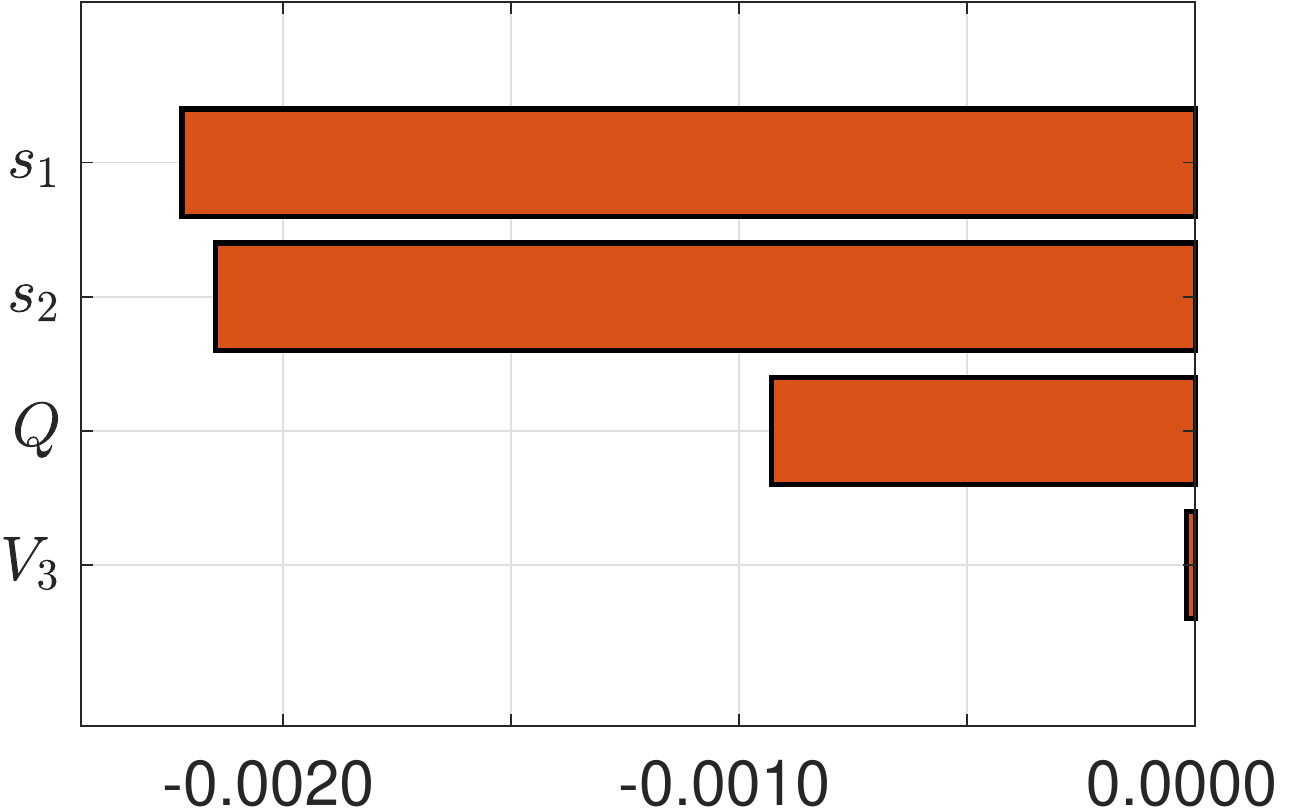}
\put(-308,112){Sensitivities for VV ECMO -- switched config. -- ADE model}
\caption{Sensitivities for VV ECMO with the switched configuration of the TPE device.  The left panel shows all of the sensitivities, and the right panel shows the four smallest sensitivities.}
\label{fig:sensVVswitchedADE}
\end{center}
\end{figure}

Figures \ref{fig:sensVAtypicalADE} and \ref{eq:VAswitchedADE} show sensitivities for the models of VA ECMO and for the typical and switched configurations of the TPE device, respectively.  Figures \ref{fig:sensVVtypicalADE} and \ref{fig:sensVVswitchedADE} show the same results but for VV ECMO with TPE.  The left panel of these figures shows results for all of the parameters, and the right panel shows results for the four parameters with the smallest sensitivity in each case.  The results reveal that for both types of ECMO and for the typical and switched configurations of the TPE device, $Q_1$ has the largest impact on the fraction of new plasma.  Further, the TPE device flow is directly related to fraction of new plasma: an increase in $Q_1$ leads to an increase in new plasma fraction.  Another feature common across ECMO types and TPE device configurations is the lack of sensitivity with respect to the ECMO compartment volume $V_3$.  Also, the fraction of new plasma is inversely related to $V_3$, meaning that a larger extracoporeal volume results in less efficient plasma exchange. 

One puzzling aspect is the sensitivity of the models to the fraction $\alpha$ of cardiac output $Q$ supported by the ECMO circuit.  For VA ECMO and both TPE device configurations, an increase in $\alpha$, corresponding to an increase in ECMO circuit flow and hence flow by the TPE device, leads to an increase in the fraction of new plasma.  This result is also seen for VV ECMO and the switched configuration of the TPE device.  Surprisingly, the opposite result is seen for VV ECMO with the typical configuration of the TPE device: an increase in ECMO circuit flow results in very small drop in fraction of new plasma. 
  
\section{Conclusions}

In this paper, we derived models for therapeutic plasma exchange done simultaneously with either veno--venous or veno--arterial extracorporeal membrane oxygentation.  These models describe fraction of new plasma in blood over time.  Two types of models were developed: algebraic delay equations and systems of delay differential equations.  We showed that our models in special cases reduce to the one described by Randerson et al.~\cite{Randerson82}, which is important because their model was validated with clinical data.  Our models extend their work by incorporating separate compartments for the native circulation and the ECMO device, and our models also account for recirculation of new plasma.  We showed that the algebraic delay equations are forward Euler discretizations of the delay differential equations, with timesteps equal to compartment transit times.  These transit times are small compared to the plasma exchange procedure, so the algebraic delay equations and systems of delay differential equation gave very similar results. Our models were applied to the scenario of port switching in the plasma exchange device, and our numerical results demonstrated that the switched configuration is less efficient in our models in exchanging plasma when the ECMO flow is very small.  Sensitivity analysis revealed the models are most sensitive to TPE device flow and least sensitive to ECMO circuit volume.

\section{Acknowledgements}

This work was supported in part by the Research Training Group in Modeling and Simulation funded by the National Science Foundation via grant RTG/DMS-1646339.

\newpage
\bibliographystyle{plain}
\bibliography{TPE_refs}

\begin{thebibliography}{10}

\bibitem{Belousova19}
Tatiana Belousova, Yi~Tong, Yu~Bai, Kimberly Klein, Hlaing Tint, and Brian
  Castillo.
\newblock Utilization of therapeutic plasma exchange for hyperbilirubinemia in
  a premature newborn on extracorporeal membrane oxygenation.
\newblock {\em Journal of Clinical Apheresis}, 34(5):615--622, 2019.

\bibitem{Chhibber12}
Vishesh Chhibber and Robert Weinstein.
\newblock Evidence-based review of therapeutic plasma exchange in neurological
  disorders.
\newblock In {\em {Seminars in Dialysis}}, volume~25, pages 132--139. Wiley
  Online Library, 2012.

\bibitem{Chong17}
Mei Chong, Alejandro~J Lopez-Magallon, Lucas Saenz, Mahesh~S Sharma, Andrew~D
  Althouse, Victor~O Morell, and Ricardo Munoz.
\newblock Use of therapeutic plasma exchange during extracorporeal life support
  in critically ill cardiac children with thrombocytopenia-associated
  multi-organ failure.
\newblock {\em Frontiers in Pediatrics}, 5:254, 2017.

\bibitem{Duan20}
Kai Duan, Bende Liu, Cesheng Li, Huajun Zhang, Ting Yu, Jieming Qu, Min Zhou,
  Li~Chen, Shengli Meng, Yong Hu, et~al.
\newblock Effectiveness of convalescent plasma therapy in severe {COVID-19}
  patients.
\newblock {\em Proceedings of the National Academy of Sciences},
  117(17):9490--9496, 2020.

\bibitem{Dyer14}
Mitchell Dyer, Matthew~D Neal, Marian~A Rollins-Raval, and Jay~S Raval.
\newblock Simultaneous extracorporeal membrane oxygenation and therapeutic
  plasma exchange procedures are tolerable in both pediatric and adult
  patients.
\newblock {\em Transfusion}, 54(4):1158--1165, 2014.

\bibitem{Jhang07}
Jeffrey Jhang, William Middlesworth, Rose Shaw, Kevin Charette, Joey Papa,
  Rashida Jefferson, Antonio~S Torloni, and Joseph Schwartz.
\newblock Therapeutic plasma exchange performed in parallel with extra
  corporeal membrane oxygenation for antibody mediated rejection after heart
  transplantation.
\newblock {\em Journal of Clinical Apheresis}, 22(6):333--338, 2007.

\bibitem{Keith20}
Philip Keith, Matthew Day, Linda Perkins, Lou Moyer, Kristi Hewitt, and Adam
  Wells.
\newblock A novel treatment approach to the novel coronavirus: an argument for
  the use of therapeutic plasma exchange for fulminant {COVID-19}.
\newblock {\em Critical Care}, 24(128), 2020.

\bibitem{Kellogg88}
Robert~M Kellogg and Jeane~P Hester.
\newblock Kinetics modeling of plasma exchange: intra-and post-plasma exchange.
\newblock {\em Journal of Clinical Apheresis}, 4(4):183--187, 1988.

\bibitem{kesici20}
Selman Kesici, Sinan Yavuz, and Benan Bayrakci.
\newblock Get rid of the bad first: Therapeutic plasma exchange with
  convalescent plasma for severe {COVID-19}.
\newblock {\em Proceedings of the National Academy of Sciences}, 2020.

\bibitem{Laverdure18}
Florent Laverdure, Laurent Masson, Guillaume Tachon, Julien Guihaire, and
  Fran{\c{c}}ois Stephan.
\newblock Connection of a renal replacement therapy or plasmapheresis device to
  the ecmo circuit.
\newblock {\em Asaio Journal}, 64(1):122--125, 2018.

\bibitem{Madore96}
Francois Madore, J~Michael Lazarus, and Hugh~R Brady.
\newblock Therapeutic plasma exchange in renal diseases.
\newblock {\em Journal of the American Society of Nephrology}, 7(3):367--386,
  1996.

\bibitem{Niemann02}
Claus~U Niemann, C~Spencer Yost, Susan Mandell, and Thomas~K Henthorn.
\newblock Evaluation of the splanchnic circulation with indocyanine green
  pharmacokinetics in liver transplant patients.
\newblock {\em Liver Transplantation}, 8(5):476--481, 2002.

\bibitem{Pavlushkov17}
Evgeny Pavlushkov, Marius Berman, and Kamen Valchanov.
\newblock Cannulation techniques for extracorporeal life support.
\newblock {\em Annals of Translational Medicine}, 5(4), 2017.

\bibitem{Puelz20}
Charles Puelz, Jonathan~L. Marinaro, Yara~A. Park, Boyce~E. Griffith,
  Charles~S. Peskin, and Jay~S. Raval.
\newblock Mathematical modeling of the impact of recirculation on exchange
  kinetics in tandem extracorporeal membrane oxygenation and therapeutic plasma
  exchange.
\newblock {\em Journal of Clinical Apheresis}, 2020.

\bibitem{Randerson82}
David~H Randerson, Matthias Blumenstein, Rupert Habersetzer, Walter Samtleben,
  Baerbel Schmidt, and Hans~J Gurland.
\newblock Mass transfer in membrane plasma exchange.
\newblock {\em Artificial Organs}, 6(1):43--49, 1982.

\bibitem{Weiss09}
Michael Weiss.
\newblock Cardiac output and systemic transit time dispersion as determinants
  of circulatory mixing time: a simulation study.
\newblock {\em Journal of Applied Physiology}, 107(2):445--449, 2009.

\bibitem{Zanella16}
Alberto Zanella, Domenico Salerno, Vittorio Scaravilli, Marco Giani, Luigi
  Castagna, Federico Magni, Eleonora Carlesso, Paolo Cadringher, Michela
  Bombino, Giacomo Grasselli, et~al.
\newblock A mathematical model of oxygenation during venovenous extracorporeal
  membrane oxygenation support.
\newblock {\em Journal of Critical Care}, 36:178--186, 2016.

\end{thebibliography}

\end{document}